\documentclass[final,onefignum,onetabnum]{siamart171218}

\usepackage{xcolor}
\usepackage{amsfonts}
\usepackage{graphicx}
\usepackage{epstopdf}
\usepackage{overpic}
\usepackage{amsmath,esint}
\usepackage{upquote,enumitem}

\ifpdf
  \DeclareGraphicsExtensions{.eps,.pdf,.png,.jpg}
\else
  \DeclareGraphicsExtensions{.eps}
\fi

\newsiamremark{remark}{Remark}
\newsiamremark{hypothesis}{Hypothesis}
\crefname{hypothesis}{Hypothesis}{Hypotheses}
\newsiamthm{claim}{Claim}

\headers{Computing with functions in the ball}{N. Boull\'e and A. Townsend}

\title{Computing with functions in the ball\thanks{Submitted to the editors \today.
\funding{The work of the first author was supported by the EPSRC Centre for Doctoral Training in Industrially Focused Mathematical Modelling (EP/L015803/1). This work is supported by National Science Foundation grant no. 1818757.}}}

\author{Nicolas Boull\'e\thanks{Mathematical Institute, University of Oxford, Oxford, OX2 6GG, UK. (\email{boulle@maths.ox.ac.uk})}
\and Alex Townsend\thanks{Department of Mathematics, Cornell University, Ithaca, NY 14853. (\email{townsend@cornell.edu})}}

\usepackage{amsopn}
\DeclareMathOperator{\diag}{diag}
 
\begin{document}

\maketitle

\begin{abstract}
A collection of algorithms in object-oriented MATLAB is described for numerically computing with smooth functions defined on the unit ball in the Chebfun software. Functions are numerically and adaptively resolved to essentially machine precision by using a three-dimensional analogue of the double Fourier sphere method to form ``ballfun" objects. Operations such as function evaluation, differentiation, integration, fast rotation by an Euler angle, and a Helmholtz solver are designed. Our algorithms are particularly efficient for vector calculus operations, and we describe how to compute the poloidal-toroidal and Helmholtz--Hodge decomposition of a vector field defined on the ball.
\end{abstract}



\section{Introduction}\label{sec_intro}
Three-dimensional spherical geometries are common in computational science and engineering, arising in weather forecasting~\cite{coiffier2011fundamentals}, geophysics~\cite{flyer2009radial,orszag1974fourier,spotz1998fast,wright2010hybrid}, hydrodynamics~\cite{hollerbach1994magnetohydrodynamic,kuang1999numerical,zhang2000magnetohydrodynamics,zhang1989convection}, and computational fluid dynamics~\cite{kerswell2005recent,serre2001three}. In each of these applications, it is routine to derive models that are continuous, even though one immediately discretizes them to compute a solution. Ballfun is a software system written in MATLAB that exploits object-oriented programming to allow users to compute with scalar- and vector-valued functions defined on the three-dimensional unit ball while being oblivious to our underlying discretizations. Ballfun is the first extension of Chebfun to three-dimensional spherical geometries~\cite{Driscoll2014} and follows the development of Spherefun~\cite{townsend2016computing} and Diskfun~\cite{wilber2017computing} for computing with functions in the sphere and the unit disk. Software systems in Dedalus~\cite{burns2019dedalus} (written in Python) and Approxfun~\cite{olver2019approxfun} (written in Julia) for computations on the ball may follow soon. Dedalus and Approxfun already have excellent functionality for computing on the 2-sphere and disks~\cite{olver2019approxfun,vasil2019tensor}.

For computations with functions defined on the unit ball, a standard approach is to employ spherical coordinates $(r,\lambda,\theta)\in[0,1]\times[-\pi,\pi]\times[0,\pi]$,
where $r$, $\lambda$, and $\theta$ denote the radial, azimuthal, and polar variables, respectively. Thus, computations on the unit ball can be conveniently related to analogous tasks involving functions defined on a cuboid, which allows for efficient algorithms based on tensor-product structure. Unfortunately, this simple coordinate transform comes with several significant disadvantages due to the artificial pole singularities introduced by the transform. 

In this paper, we employ a technique known as the double Fourier sphere (DFS) method~\cite{fornberg1997comparison,merilees1973pseudospectral,orszag1974fourier} in conjunction with tensor-product expansions of functions. More precisely, we use a three-dimensional analogue of the DFS method that extends ideas from the disk and sphere~\cite{townsend2016computing,wilber2017computing} (implemented in Spherefun and Diskfun, which are part of Chebfun), while preserving the additional structure that is present in the 3D ball (see~\cref{def:BMC3}). The DFS method alleviates some of the computational difficulties with spherical coordinates while having approximants that have an underlying tensor-product structure for efficient algorithms and FFT-based fast transforms. We use DFS approximants to develop a collection of algorithms for performing everyday computational tasks on scalar- and vector-valued functions defined on the unit ball and thus provide a convenient computational environment for scientific explorations.    

Our algorithms are designed, whenever possible, to compute each operation on a function to essentially machine precision by using data-driven techniques from approximation theory.  Our codes are also designed to have no required user-defined algorithmic parameters and to be as intuitive as possible for MATLAB users. For example, \texttt{sum(v)} returns the sum of the entries of a vector \texttt{v} in MATLAB, while \texttt{sum3(f)} returns the integral of a function \texttt{f} over the unit ball. Moreover, \texttt{v.*w} performs entry-by-entry vector multiplication, while \texttt{f.*g} returns a function representing the multiplication of \texttt{f} and \texttt{g} in Ballfun.  During the operation $\texttt{f.*g}$, our algorithm automatically selects the discretization of the output so that the result is approximated to essentially machine precision.  We repeat this idea in the one hundred or so Ballfun commands by constantly expanding and pruning underlying discretizations to represent functions as efficiently as possible.  

There are several existing approaches for computing with functions on the unit ball and we seriously considered two other approaches: 
\begin{description}[wide=\parindent]
\item[Spherical harmonic expansions:]  Spherical harmonic expansions of a function are given by $f(r,\lambda,\theta) =  \sum_{\ell=0}^\infty \sum_{m=-\ell}^\ell f_\ell^m r^\ell Y_{\ell}^m(\lambda, \theta)$, where $Y_{\ell}^m$ is a spherical harmonic. These expansions can be thought of as the ball analogue of trigonometric expansions for periodic functions. When truncated, they provide essentially uniform resolution of a function over the ball. They have major applications in geophysics~\cite{lowrie2007fundamentals} and the numerical solution of separable elliptic equations. 

\item[Orthogonal polynomials on the ball:] Given an appropriate weight function on the ball, one can derive various families of orthogonal polynomials that are built from ultraspherical polynomials~\cite[Sec.~5.1]{dunkl2014orthogonal}. Expanding functions in any one of these bases provides excellent resolution properties, along with fast evaluation, differentiation, and integration of the expansions. Unlike spherical harmonics expansions, they are rarely employed in practice. 
\end{description}

We require a representation for functions on the ball that can be adaptively computed, as we would like to achieve an accuracy close to machine precision. While there are optimal-complexity spherical harmonic transforms~\cite{slevinsky2017fast}, it is highly desirable to have the most computationally efficient fast transform associated with an expansion. The DFS method offers a simple and computationally efficient fast transform based on the FFT (see~\cref{sec_compute_coeffs}). Unlike spherical harmonic and orthogonal polynomial expansions, the DFS method does not guarantee that an expansion is infinitely differentiable on the ball, even when the original function is. For this reason, our algorithms must strictly preserve a structure in the DFS expansions to ensure that they represent a smooth function on the ball (see~\cref{def:BMC3}). 

Using DFS approximants, we develop a variety of algorithmic tools to provide a convenient computational environment for integrating, differentiating, and solving partial differential equations (see~\cref{sec_Helmholtz}), as well as representing vector-valued functions. This allows us to develop a set of algorithms for performing vector calculus (see~\cref{sec_vector}), including computing the Helmholtz--Hodge and poloidal-toroidal decomposition. 

The paper is organized as follows.  We briefly introduce the software that accompanies this paper in~\cref{sec:software}. Then, in~\cref{sec_constructor}, we explain the methods used to discretize smooth functions on the ball. Next, in \cref{sec_algo}, we discuss some of the operations implemented in the software such as integration, differentiation, and a fast rotation algorithm. Following this, in \cref{sec_Helmholtz}, we describe a fast and spectrally accurate Helmholtz solver for solving equations with Dirichlet or Neumann boundary conditions. Finally, \cref{sec_vector} consists of a description of the vector calculus algorithms, including the poloidal-toroidal and Helmholtz--Hodge decompositions.

\subsection{Software}\label{sec:software}
Ballfun is part of Chebfun~\cite{Driscoll2014}, which is a software system for computing with functions and solving differential equations on an interval~\cite{battles2004extension}, rectangle~\cite{townsend2013extension}, cuboid~\cite{hashemi2017chebfun}, disk~\cite{wilber2017computing}, and the surface of a sphere~\cite{townsend2016computing}.  Accompanying this paper is the publicly available MATLAB code in Chebfun~\cite{Driscoll2014} with two new classes called \texttt{ballfun} and \texttt{ballfunv}. We encourage the reader to explore this paper with the latest version of Chebfun downloaded and ready for interactive exploration. On the Chebfun website, we provide documentation in the form of a chapter of the Chebfun Guide~\cite{Driscoll2014} as well as several examples.\footnote{The Ballfun examples are available at \url{http://www.chebfun.org/examples/sphere/}.} Functions on the ball can be easily constructed in the software by calling the appropriate command. For instance,
\texttt{f = ballfun(@(x,y,z) sin(cos(y)))} defines the function $f(x,y,z) = \sin(\cos y)$.
Underneath, Ballfun adaptively resolves the function to machine precision and represents it using the DFS method. For example,
\begin{verbatim}
  f = ballfun(@(x,y,z) sin(cos(y))) % ballfun representing sin(cos(y))
  ballfun object:
  	    domain           r    lambda    theta
      unit ball        21      45       41
\end{verbatim}
where $21$, $45$, and $41$ are discretization parameters that Ballfun automatically determined necessary to resolve $f$ to machine precision. The Ballfun software is highly adaptive and automatically truncates the expansion to resolve functions on the ball to machine precision after each operation. After its construction (see~\cref{sec_constructor}), a function can be manipulated and analyzed through the nearly one hundred operations implemented in the package (see~\cref{tab:scalar} and~\cref{tab:vector}). 

\begin{table}[htbp]\label{tab:scalar}  
\centering
\caption{A selection of Ballfun commands for scalar-valued functions.}
\begin{tabular}{cc}
Ballfun command & Operation \\
\hline \\[-8pt]
\texttt{+}, \texttt{-}, \texttt{.*}, \texttt{./} & basic arithmetic \\ 
\texttt{coeffs2vals}, \texttt{vals2coeffs} & fast transforms\\
 \texttt{feval} & pointwise evaluation \\ 
 \texttt{sum}, \texttt{sum2}, \texttt{sum3} & integration \\
 \texttt{diff} & differentiation \\
 \texttt{rotate} & rotate using Euler angles \\
 \texttt{helmholtz} & helmholtz solver \\
\end{tabular} 
\end{table} 

\begin{table}[htbp]\label{tab:vector}  
\centering
\caption{A selection of Ballfun commands for vector-valued functions.}
\begin{tabular}{cc}
Ballfun command & Operation \\
\hline \\[-8pt]
\texttt{cross} & cross-product \\
\texttt{dot} & dot-product \\ 
\texttt{feval} & pointwise evaluation \\ 
\texttt{curl} & curl \\ 
\texttt{divergence} & divergence \\
\texttt{PTdecomposition} & poloidal-toroidal decomposition \\
\texttt{HelmholtzDecomposition} & Helmholtz--Hodge decomposition
\end{tabular} 
\end{table} 

\section{The Ballfun constructor}\label{sec_constructor}
In this section, we explain how smooth functions, expressed in Cartesian or spherical coordinates, are discretized and constructed in our software. Smooth functions on the ball expressed in the spherical coordinate system $(r,\lambda,\theta)\in[0,1]\times[-\pi,\pi]\times[0,\pi]$ can potentially introduce artificial boundaries at the origin or poles, as well as the loss of periodicity in the polar variable $\theta$. To overcome this issue, we first sample functions on a tensor-product grid in spherical coordinates. Then, we compute a Chebyshev--Fourier--Fourier expansion that interpolates the samples, using the ball analogue for the double Fourier sphere (DFS) method.

\subsection{The double Fourier sphere method in the ball} 
\label{sec_DFS}
The DFS method for the sphere was originally proposed by Merilees~\cite{merilees1973pseudospectral} and is used to construct spherefun objects in Chebfun. It naturally extends to the 3D settings and maps a function defined on a ball onto a 3D cuboid so that the origin and poles of the ball are not treated as artificial boundaries and the polar variable can be represented in a Fourier series~\cite{boyd1978choice,fornberg1995pseudospectral,orszag1974fourier,yee1980studies}. The method can also be applied to disks, cylinders, and ellipsoids~\cite{townsend2016computing,wilber2017computing}.

The ball analogue of the DFS method is obtained by constructing a Chebyshev--Fourier--Fourier expansion of a function defined on $[-1,1]\times [-\pi,\pi]\times [-\pi,\pi]$ instead of $[0,1]\times [-\pi,\pi]\times [0,\pi]$. A continuous function $f(x,y,z)$ on the ball is first written in spherical coordinates as
\[f(r,\lambda,\theta)=f(r\cos\lambda\sin\theta,r\cos\lambda\cos\theta,r\cos\theta),\quad (r,\lambda,\theta)\in [0,1]\times [-\pi,\pi]\times [0,\pi].\]
The function $f(r,\lambda,\theta)$ is not periodic in $\theta$. Under the DFS mapping, it is recovered by ``doubling up'' the polar variable to $[-\pi,\pi]$ in the sense that $f$ is sampled twice. The radial variable is also doubled to remove the artificial boundary at $r=0$. Using this ideas, we extend the function $f$ to a new function $\tilde{f}$, defined on $[-1,1]\times [-\pi,\pi]\times [-\pi,\pi]$. The function $\tilde{f}$ can be expressed as
\begin{equation}
\label{eq_BMC}
\tilde{f}(r,\lambda,\theta) = \left\{
    \begin{array}{llllllll}
        g(r,\lambda+\pi,\theta), & (r,\lambda,\theta)\in [0,1]\times [-\pi,0]\times [0,\pi],\\
        h(r,\lambda,\theta), & (r,\lambda,\theta)\in [0,1]\times [0,\pi]\times [0,\pi],
        \\ \\
        g(-r,\lambda,\pi-\theta), & (r,\lambda,\theta)\in [-1,0]\times [0,\pi]\times [0,\pi],\\
        h(-r,\lambda+\pi,\pi-\theta), & (r,\lambda,\theta)\in [-1,0]\times [-\pi,0]\times [0,\pi],\\ \\
        h(r,\lambda+\pi,-\theta), & (r,\lambda,\theta)\in [0,1]\times [-\pi,0]\times [-\pi,0],\\
        g(r,\lambda,-\theta), & (r,\lambda,\theta)\in [0,1]\times [0,\pi]\times [-\pi,0],
        \\ \\
        h(-r,\lambda,\pi+\theta), & (r,\lambda,\theta)\in [-1,0]\times [0,\pi]\times [-\pi,0],\\
        g(-r,\lambda+\pi,\pi+\theta), & (r,\lambda,\theta)\in [-1,0]\times [-\pi,0]\times [-\pi,0],
    \end{array}
\right.
\end{equation}
where
\[g(r,\lambda,\theta) = f(r,\lambda-\pi,\theta),\quad h(r,\lambda,\theta) = f(r,\lambda,\theta), \quad (r,\lambda,\theta)\in[0,1]\times[0,\pi]\times[0,\pi].\]
Functions that satisfy~\cref{eq_BMC} are said to be block-mirror-centrosymmetric (BMC)~\cite{wilber2017computing}. A more intuitive description is given by the visualization
\begin{equation}
\label{eq_flip}
\tilde{f} = \Bigg[
\begin{bmatrix}
g & h\\
\texttt{flip1}(\texttt{flip3}(h)) & \texttt{flip1}(\texttt{flip3}(g))
\end{bmatrix};
\begin{bmatrix}
\texttt{flip3}(h) & \texttt{flip3}(g)\\
\texttt{flip1}(g) & \texttt{flip1}(h)
\end{bmatrix}\Bigg],
\end{equation}
where \texttt{flip1} (resp. \texttt{flip3}) refers to the MATLAB command that reverses the order of the first (resp. third) component of a tensor.

In addition to satisfying the BMC structure, $\tilde{f}$ must be constant at $r=0$ as well as $\theta = 0$ and $\theta=\pi$, corresponding to the origin and the poles. We call these functions BMC-III functions. (BMC-I and BMC-II functions are defined in~\cite{townsend2016computing} and~\cite{wilber2017computing}, respectively.)

\begin{definition}[BMC-III function]\label{def:BMC3}
A function $\tilde{f}:[-1,1]\times [-\pi,\pi]\times [-\pi,\pi]\rightarrow\mathbb{C}$ is a BMC-III (Type-III BMC) function if it is a BMC function, $\tilde{f}(0,\cdot,\cdot)=\alpha$, and, for any $r\in[0,1]$, $\tilde{f}(r,\cdot,0)=\beta(r)$ and $\tilde{f}(r,\cdot,\pi)=\gamma(r)$, where $\beta$ and $\gamma$ only depend on $r$ such that $\beta(0)=\gamma(0)=\alpha$, for some constant $\alpha$.
\end{definition}

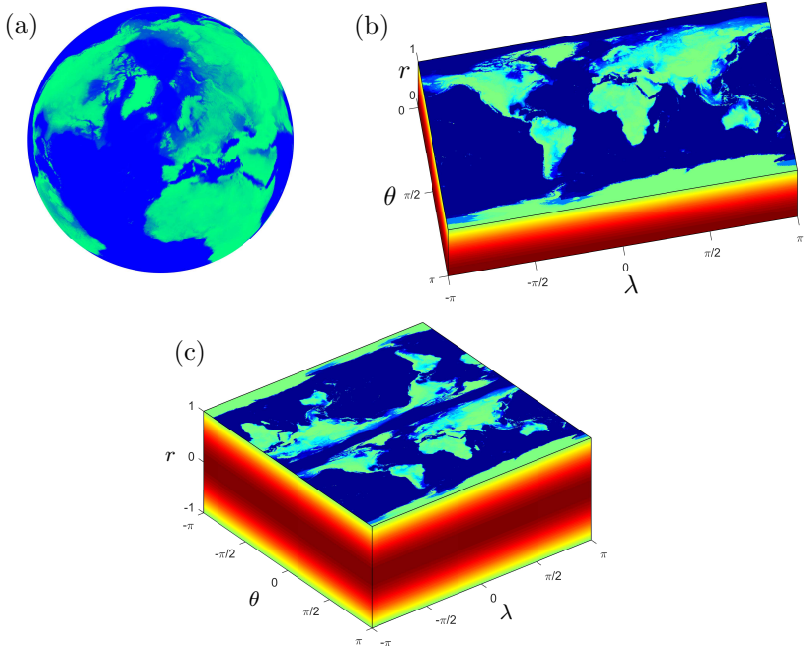
\begin{figure}[htbp]
\centering
\begin{minipage}{.4\textwidth} 
\begin{overpic}[width=\textwidth,trim={90 70 90 65},clip]{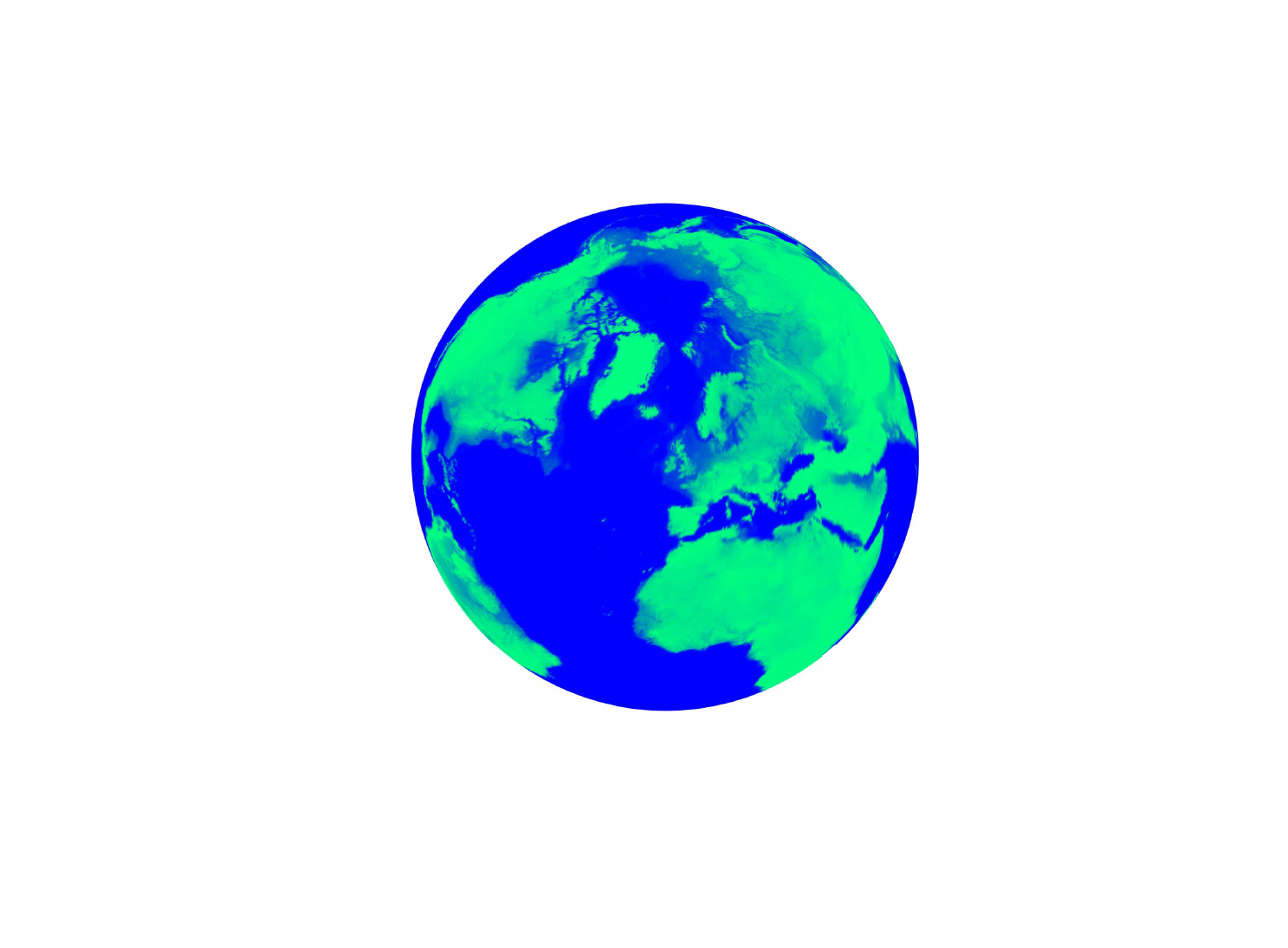}
\put(13,65){(a)}
\end{overpic}
\end{minipage}
\begin{minipage}{0.5\textwidth} 
\begin{overpic}[width=\textwidth,trim={0 10 0 10},clip]{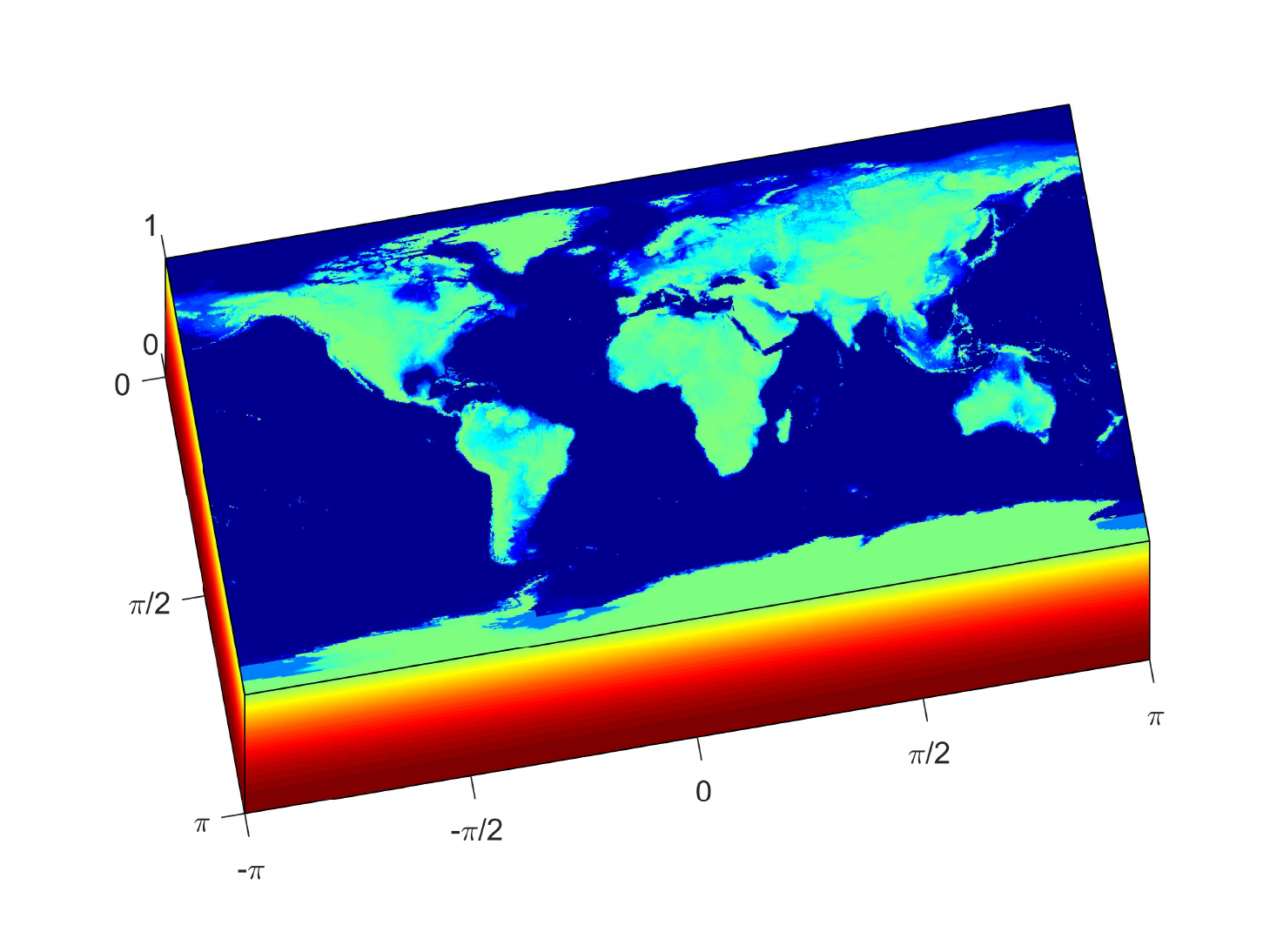}
\put(0,58){(b)}
\put(9,49){$r$}
\put(6,23){$\theta$}
\put(55,5){$\lambda$}
\end{overpic} 
\end{minipage}\\
\begin{minipage}{.45\textwidth} 
\begin{overpic}[width=\textwidth,trim={30 20 20 20},clip]{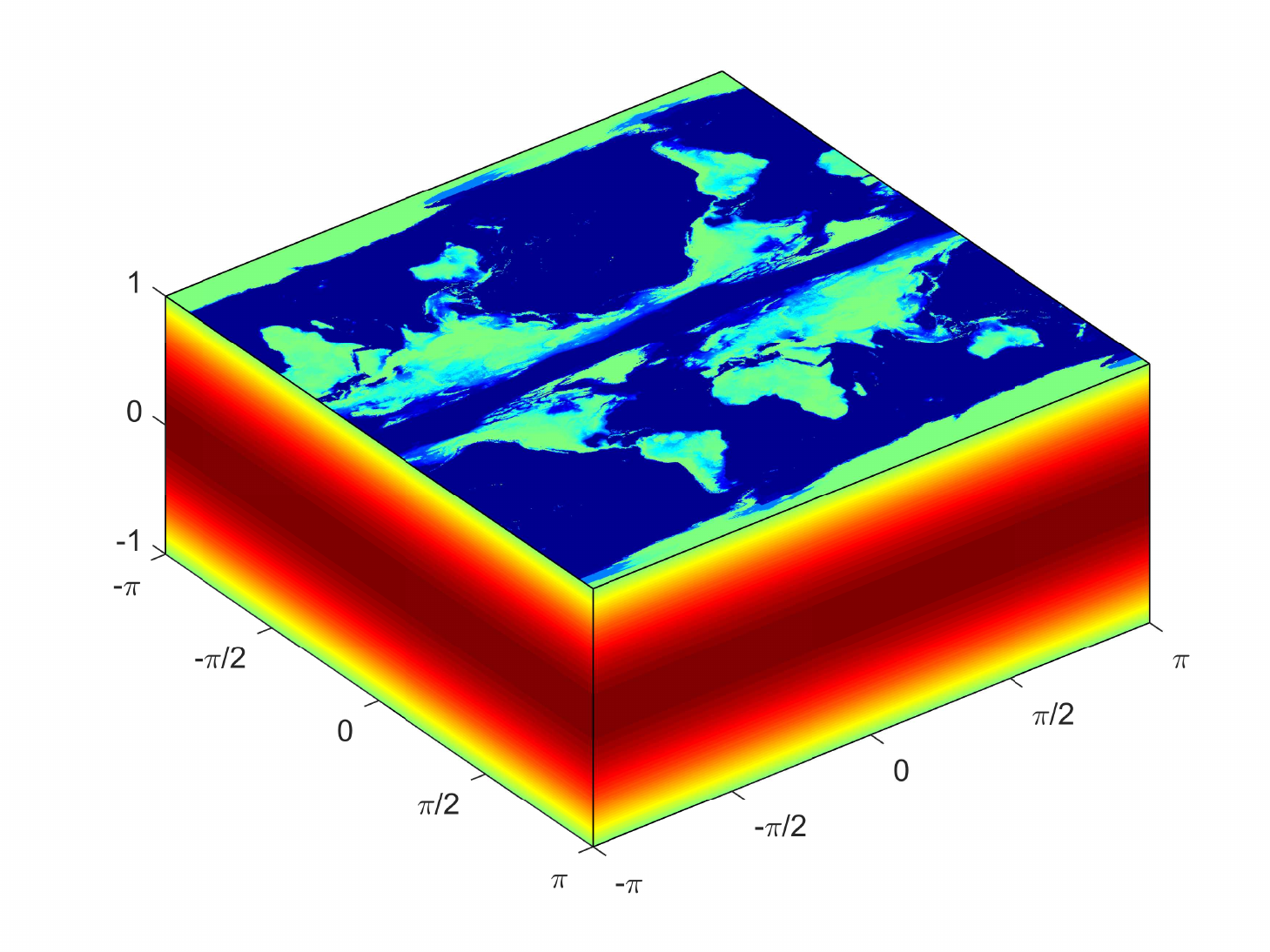}
\put(0,65){(c)}
\put(-2,42){{\footnotesize $r$}}
\put(17,9){{\footnotesize $\theta$}}
\put(74,6){{\footnotesize $\lambda$}}
\end{overpic} 
\end{minipage}
\label{fig_DFS}
\caption{The DFS method applied to the globe. (a) The solid earth including the land masses. (b) The projection of the land masses using spherical coordinates. (c) Land masses after applying the DFS method. This is a BMC-III function that is periodic in $\lambda$ and $\theta$ and defined over $(r,\lambda,\theta)\in[-1,1]\times [-\pi,\pi]\times [-\pi,\pi]$.}
\end{figure}

\Cref{fig_DFS} shows the DFS method applied to the earth and the type-III BMC structure.

There are two salient features of the DFS method that make it attractive for developing a package for computing with functions on the ball. First, tensor product expansions of Fourier and Chebyshev bases can be used to represent $\tilde{f}$. If $f(x,y,z)$ is a function in Cartesian coordinates on the ball, then after applying the DFS method, we have a function $\tilde{f}(r,\lambda,\theta)$ defined on the cuboid $(r,\lambda,\theta)\in[-1,1]\times[-\pi,\pi]\times[-\pi,\pi]$ that can be approximated as
\begin{equation}
\label{eq:expansion}
\tilde{f}(r,\lambda,\theta) \approx \sum_{i=0}^{n-1}\sum_{j=-n/2}^{n/2-1}\sum_{k=-n/2}^{n/2-1} \alpha_{ijk} T_{i}(r)e^{\mathbf{i}j\lambda}e^{\mathbf{i}k\theta},
\end{equation}
where $(r,\lambda,\theta)$ are spherical coordinates, $T_i$ denotes the Chebyshev polynomial of the first kind of degree $i$, and $n$ is an even integer that is determined by the adaptive procedure described in \cref{sec_determination_size}.\footnote{For simplicity, we have chosen the same number of terms in each sum. Ballfun uses a generalization of~\cref{eq:expansion} with a different number of terms in each sum.} This representation allows us to use fast transforms as well as 1D and 2D algorithms for Chebyshev and trigonometric expansions. The second feature is that the DFS mapping of a function leads to a BMC structure (see~\cref{fig_DFS}) that, if preserved, ensures smoothness of the solution throughout the entire domain~\cite{townsend2016computing,wilber2017computing}. The BMC symmetry is imposed exactly by evaluating the function on $(r,\lambda,\theta)\in[0,1]\times[0,\pi]\times[0,\pi]$ and extending it to $[-1,1]\times [-\pi,\pi]\times [-\pi,\pi]$ using \cref{eq_flip}. This ensures that the resulting function is smooth on the ball, i.e., at least continuous and differentiable.  There are representations of functions on the ball that preserve full regularity~\cite{vasil2019tensor}, however there are less appropriate in our settings since they do not allow efficient FFT-based transforms.

\subsection{Computing the Chebyshev--Fourier--Fourier coefficients} \label{sec_compute_coeffs}
Once the BMC-III function $\tilde{f}$ is found, it is approximated by a truncated Chebyshev--Fourier--Fourier (CFF) series~\cite{mason2002chebyshev,trefethen2000spectral,trefethen2013approximation}. For some even integer $n$, the CFF coefficients are stored as an $n\times n\times n$ tensor and the entries are computed in $\mathcal{O}(n^3\log n)$ operations, as follows:

\begin{enumerate}[wide=\parindent]
\item The function $f$ is evaluated over $[0,1]\times [-\pi,\pi]\times [0,\pi]$ at the tensor-product grid:
\begin{equation} 
\left(\!\cos\!\left(\!\frac{\left(\frac{n}{2}-1-i\right)\pi}{n-1}\!\right)\!,\frac{2j\pi}{n},\frac{2k\pi}{n}\!\right)\!, \!\quad 0\leq i\leq \frac{n}{2}-1, \!\quad - \frac{n}{2}\leq j\leq \frac{n}{2}-1, \!\quad 0\leq k\leq  \frac{n}{2}.
\label{eq:grid} 
\end{equation} 
\item The samples of $f$ are doubled-up (see~\cref{eq_flip}). This extends them to be samples of $\tilde{f}$ on $[-1,1]\times [-\pi,\pi]\times [-\pi,\pi]$ at the tensor-product grid:
\[
\left(\cos\left(\frac{i\pi}{n-1}\right),\frac{2j\pi}{n},\frac{2k\pi}{n}\right),\quad 0\leq i\leq n-1, \quad -\frac{n}{2}\leq j\leq \frac{n}{2}-1, \quad -\frac{n}{2}\leq k\leq \frac{n}{2}-1
\]
without any additional evaluations of $f$. 
\item The CFF coefficients are computed using the discrete Chebyshev tranform (DCT) \cite{gentleman1972implementing,gentleman1972implementing2,mason2002chebyshev} and the FFT~\cite{cooley1965algorithm}.
\end{enumerate}

There is also the inverse procedure, which evaluates $f$ at the grid in~\eqref{eq:grid} in $\mathcal{O}(n^3\log n)$ operations. This operation is particularly important in our plotting commands and is achieved by reversing steps 1-3, using the inverse DCT and FFT.

\subsection{Determination of the discretization size} \label{sec_determination_size}
To construct a ballfun object to represent a given function $f$, we first sample $\tilde{f}$ at a $17\times 17\times 17$ CFF grid and compute the corresponding CFF coefficients (see~\cref{sec_compute_coeffs}). We then successively increase the grid size independently in each variable from $17$ to $33$ to $65$, and so on, until we deem the function to be resolved in each variable. We use these samples to compute the CFF coefficients $A = (\alpha_{ijk})$ corresponding to an $m\times n\times p$, where $m,n,p=17, 33, 65,\ldots,$ and then gauge the resolution in each variable by creating vectors of the absolute maximum of the coefficients along each variable, i.e.,
\[
\text{Cols}_i = \max_{j,k} \left|\alpha_{ijk}\right|,\quad \text{Rows}_j = \max_{i,k}\left|\alpha_{ijk}\right|,\quad \text{Tubes}_k = \max_{i,j}\left|\alpha_{ijk}\right|.
\]
One can now inspect these vectors to identify whether or not the function is resolved to machine precision in each variable, relative to the magnitude of $f$ on $[0,1]\times [-\pi,\pi]\times [0,\pi]$~\cite{aurentz2017chopping,hashemi2017chebfun,wright2015extension}.  One can identify a near-optimal discretization size in that variable by recording the last entry in each vector above machine precision~\cite{aurentz2017chopping}.

\begin{figure}[htbp]
\centering
\begin{overpic}[width=\textwidth,trim={100 20 100 30},clip]{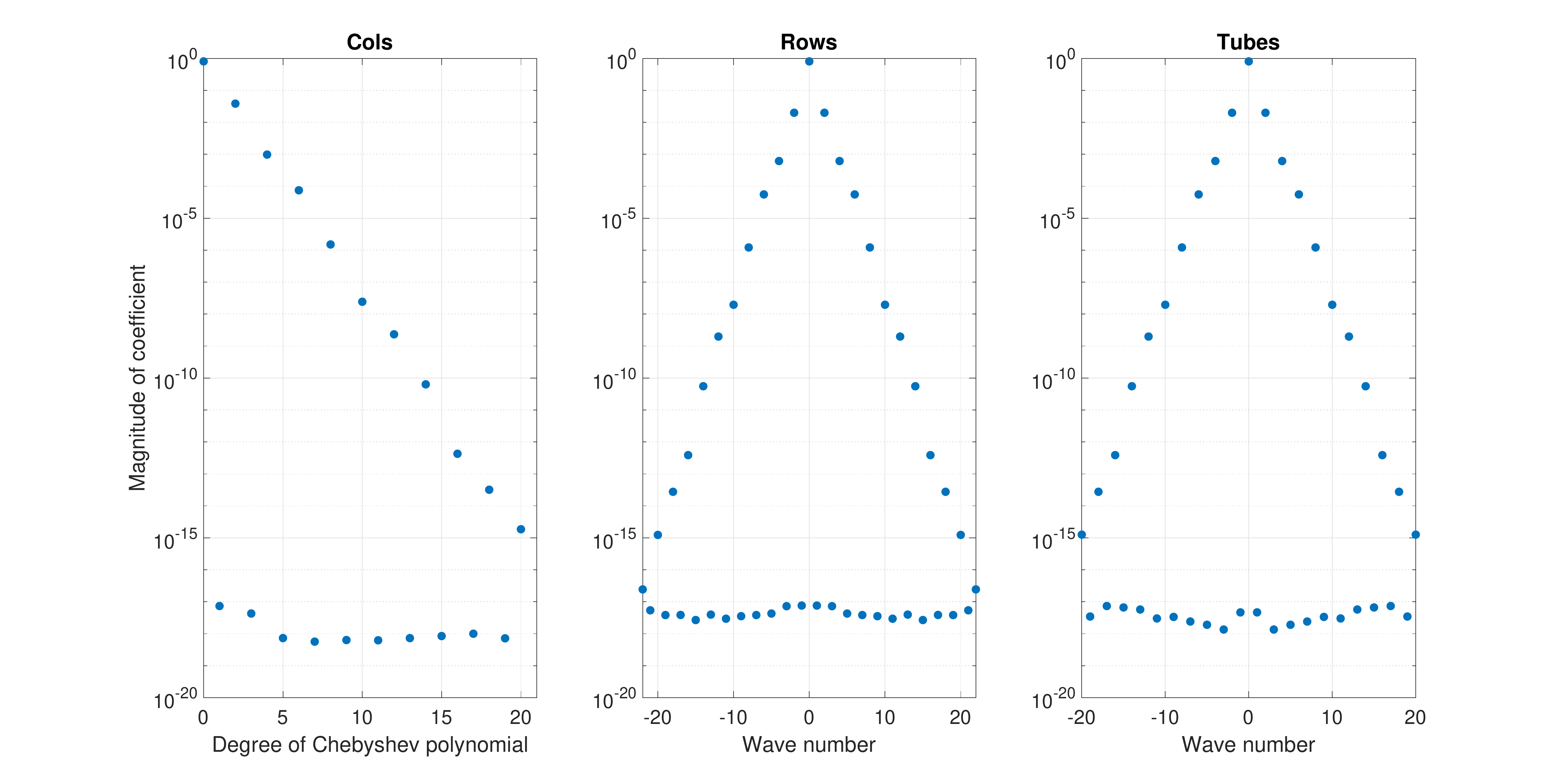}
\end{overpic} 
\label{fig_truncation}
\caption{The absolute maximum Chebyshev and Fourier coefficients in the radial, azimuthal, and polar variables of $f(x,y,z)=\sin(\cos y)$. The Ballfun constructor selected a discretization size of $21\times 45\times 41$ to represent $f$ to essentially machine precision over the ball. One can visually see that the function is likely to be resolved as the entries of \text{Cols}, \text{Rows}, and \text{Tubes} decay to machine precision.}
\end{figure}

The constructor typically terminates when it encounters vectors \text{Cols}, \text{Rows}, and \text{Tubes} as shown in~\cref{fig_truncation} for $f(x,y,z) = \sin(\cos y)$. In particular, for $f(x,y,z) = \sin(\cos y)$, the Ballfun constructor selected a CFF series of size $21\times 45\times 41$ to represent $f$ to essentially machine precision over the ball. Once the function $\tilde{f}$ is represented in a CFF expansion, the approximant is stored as a \texttt{ballfun} object, ready for further computation. For the rest of this paper, we will assume that the functions are represented by an $n\times n \times n$ tensor (though our software can deal with rectangular discretizations).

\section{Algorithms for numerical computations with functions on the ball}\label{sec_algo}
Once a \texttt{ballfun} object is computed, there are many operations that can be performed on it. In fact, many of the operations can be decomposed into a sequence of 1D operations, which are particularly efficient for approximants of the form~\cref{eq:expansion}.  This includes pointwise evaluation (see~\cref{sec:eval}), integration (see~\cref{sec:integration}), differentiation (see~\cref{sec:diff}), and fast rotation (see~\cref{sec:rotate}). 

\subsection{Pointwise evaluation}\label{sec:eval} 
The evaluation of a function $f(r,\lambda,\theta)$ at a point $(r^*,\lambda^*,\theta^*)$ in the ball can be computed in $\mathcal{O}(n^3)$ operations. It follows from the CFF approximation of $\tilde{f}$ in~\cref{eq:expansion} that one has 
\[
\begin{aligned} 
\tilde{f}(r^*,\lambda^*,\theta^*) &= \sum_{i=0}^{n-1}\sum_{j=-n/2}^{n/2-1}\sum_{k=-n/2}^{n/2-1} \alpha_{ijk} T_{i}(r^*)e^{\mathbf{i}j\lambda^*}e^{\mathbf{i}k\theta^*} \\
&=\sum_{k=-n/2}^{n/2-1} \left( \sum_{j=-n/2}^{n/2-1}\left(\sum_{i=0}^{n-1}\alpha_{ijk} T_{i}(r^*)\right) e^{\mathbf{i}j\lambda^*} \right) e^{\mathbf{i}k\theta^*} .
\end{aligned} 
\]
Therefore, $\tilde{f}(r^*,\lambda^*,\theta^*)$ can be computed by first evaluating $\sum_{i=0}^{n-1}\alpha_{ijk} T_{i}(r^*)$ using Clenshaw's algorithm~\cite[Chap.~19]{trefethen2013approximation}, which returns an $n\times n$ matrix of values. Then, one can compute the summand over the $j$ index using Horner's scheme~\cite{wright2015extension}, which returns an $n\times 1$ vector, before finally computing the summand over the $k$ index using Horner's scheme.  This algorithm is implemented in the \texttt{feval} command and returns a scalar for $\tilde{f}(r^*,\lambda^*,\theta^*)$. It is also possible to evaluate functions in Cartesian coordinates and Ballfun does a change of variables to spherical coordinates in this case.

\subsection{Integration}\label{sec:integration} 
The triple definite integral of a function $\tilde{f}(r,\lambda,\theta)$ on the unit ball can be written as follows:
\begin{align*}
\int_{B(0,1)}\tilde{f}(r,\lambda,\theta)\,dV&=\int_0^1\int_0^\pi\int_{-\pi}^\pi \tilde{f}(r,\lambda,\theta)r^2\sin\theta \,d\lambda \,d\theta \,dr\\
&= \sum_{i=0}^{n-1}\sum_{j=-n/2}^{n/2-1}\sum_{k=-n/2}^{n/2-1}\alpha_{ijk}\int_0^1r^2T_i(r)\,dr\int_{-\pi}^\pi e^{\mathbf{i}j\lambda\,}d\lambda\int_0^\pi\!\sin\theta e^{\mathbf{i}k\theta}\,d\theta,\\
&= 2\pi\sum_{i=0}^{n-1}\sum_{k=-n/2}^{n/2-1}\alpha_{i0k} \left(\int_0^1r^2T_i(r)\,dr \right) \left( \int_0^\pi\sin\theta e^{\mathbf{i}k\theta}\,d\theta\right)\\
& = 2\pi\sum_{i=0}^{n-1}\sum_{k=-n/2}^{n/2-1}\alpha_{i0k}\nu_i \frac{1+e^{i\pi k}}{1-k^2}, \quad \nu_i = \begin{cases}\frac{2i^2-6}{i^4-10i^2+9}, & i \text{ even}, \\ 0, & i \text{ odd}.\end{cases} 
\end{align*}
Here, the last equality follows by calculating the integrals in $r$ explicitly (see, for example,~\cite[(4.3)]{townsend2016computing}).  Therefore, the integral of $\tilde{f}$ reduces to a basic task in linear algebra and can be computed in $\mathcal{O}(n^2)$ operations. This is implemented in Ballfun in the  \texttt{sum3} command. For example, the function $f(x,y,z)=x^2$ has an integral of $4\pi/15$ over the ball and can be computed in Ballfun by
\begin{verbatim}
  f = ballfun(@(x,y,z) x.^2); % ballfun representing x^2
  sum3(f)                     % Integrate over the ball 
  ans =
      0.837758040957278
\end{verbatim}
The error is computed as \texttt{abs(sum3(f) - 4*pi/15)} and is given by $1.102\times 10^{-16}$. 

\subsection{Differentiation}\label{sec:diff}
Differentiation of a function on the ball with respect to spherical coordinates in $r$, $\lambda$, and $\theta$ may introduce singularities at the poles and origin. For instance, consider the smooth function $f(r,\lambda,\theta) = r\cos\theta$. The derivative of $f$ with respect to $\theta$ is $-r\sin\theta$, which is not smooth at the poles. However, we are interested in computing derivatives that arise in vector calculus, such as the gradient, the divergence, the curl, or the Laplacian. All these operations can be expressed as partial derivatives in the Cartesian coordinates system.  Therefore, our default is to allow for ballfun objects to be differentiated in the Cartesian coordinate system. 

We follow the same approach as Spherefun~\cite{townsend2016computing} and express the partial derivatives in $x$, $y$, and $z$ in terms of the spherical coordinates $r$, $\lambda$, and $\theta$ as follows:
\begin{align}
\frac{\partial}{\partial x}&=\cos\lambda\sin\theta\frac{\partial}{\partial r}-\frac{\sin\lambda}{r\sin\theta}\frac{\partial}{\partial\lambda}+\frac{\cos\lambda\cos\theta}{r}\frac{\partial}{\partial\theta},\label{eq_diff_1}\\
\frac{\partial}{\partial y}&=\sin\lambda\sin\theta\frac{\partial}{\partial r}+\frac{\cos\lambda}{r\sin\theta}\frac{\partial}{\partial\lambda}+\frac{\sin\lambda\cos\theta}{r}\frac{\partial}{\partial\theta},\label{eq_diff_2}\\
\frac{\partial}{\partial z}&=\cos\theta\frac{\partial}{\partial r}+\frac{\sin\theta}{r}\frac{\partial}{\partial\theta}.\label{eq_diff_3}
\end{align}
Then, \cref{eq_diff_1}, \cref{eq_diff_2}, and \cref{eq_diff_3} involve $\mathcal{O}(n^3)$ operations on the tensor of CFF coefficients representing $\tilde{f}$. For example, the derivative of $\tilde{f}$ with respect to $\lambda$ can be expressed as
\[
\frac{\partial\tilde{f}}{\partial\lambda} = \sum_{i=0}^{n-1}\sum_{j=-n/2}^{n/2-1}\sum_{k=-n/2}^{n/2-1} \alpha_{ijk} \mathbf{i}jT_{i}(r)e^{\mathbf{i}j\lambda}e^{\mathbf{i}k\theta}.
\]

Multiplications and divisions by $\sin\lambda$, $\cos\lambda$, $\sin\theta$, and $\cos\theta$ in~\cref{eq_diff_1,eq_diff_2,eq_diff_3} are computed by multiplying the tensor of CFF coefficients $A=(\alpha_{ijk})$ by the corresponding matrices of linear operators, expressed in the Fourier basis. For example, we write $\tilde{f}(r,\lambda,\theta)/\sin\theta \approx \sum_{i=0}^{n-1}\sum_{j=-n/2}^{n/2-1}\sum_{k=-n/2}^{n/2-1} \beta_{ijk}\mathbf{i}jT_{i}(r)e^{\mathbf{i}j\lambda}e^{\mathbf{i}k\theta}$, where $B = (\beta_{ijk})$ satisfies 
\[
B(:,j,:) = A(:,j,:) M_{\sin}^{-\top},
\quad -\frac{n}{2}\leq j\leq \frac{n}{2}-1,
\quad M_{\sin}=\frac{\mathbf{i}}{2}
\begin{bmatrix}
 0 &      1 &      &  \\
-1 & \ddots &\ddots&  \\
   &  \ddots&\ddots& 1\\
   &        &   -1 & 0
\end{bmatrix}.
\]
Here, $M_{\sin}$ is the matrix of multiplication by $\sin \theta$ in the Fourier basis. It is nonsingular if we choose $n$ to be even (in this case the eigenvalues are $\cos(\pi l/(n+1))$, $1\leq l\leq n$). Moreover, we write  $\tilde{f}(r,\lambda,\theta)/r \approx \sum_{i=0}^{n-1}\sum_{j=-n/2}^{n/2-1}\sum_{k=-n/2}^{n/2-1} \beta_{ijk}\mathbf{i}jT_{i}(r)e^{\mathbf{i}j\lambda}e^{\mathbf{i}k\theta}$, where $B = (\beta_{ijk})$ satisfies 
\[
B(:,j,:) = M_r^{-1} A(:,j,:),
\quad -\frac{n}{2}\leq j\leq \frac{n}{2}-1,
\quad M_{r}=\
\begin{bmatrix}
0 &\frac{1}{2}&            &            & 			  &\\
1 &    		  0&\frac{1}{2} &            &  		  &\\
  &\frac{1}{2}&     \ddots &\ddots      & 			  &\\
  &           &\ddots		&\ddots		 &\frac{1}{2} &\\
  &			  & 			&\frac{1}{2} &0 	      &\frac{1}{2}\\
  &			  & 			&            &\frac{1}{2} & 0

\end{bmatrix}.
\]
Here, $M_r$ stands for the matrix of multiplication by $r$ in the Chebyshev basis. This matrix is invertible for even $n$ since its determinant is equal to $-(-1/4)^{n/2-1}/2$.

Working directly on coefficients allows us to circumvent potential singularity issues at $r=0$ and the poles, while the standard technique uses a ``shifted grid" procedure in the physical space~\cite{cheong2000application,fornberg1995pseudospectral,heinrichs2004spectral}. This procedure shifts the grid of sampled points in the latitude and radial directions, which avoids evaluation at the poles and the origin but can be numerically inaccurate near these points.

These operations are implemented in Ballfun in the \texttt{diff} command. For example, the derivative of $f(x,y,z) = \cos(xy)$ with respect to $x$ is also represented as a ballfun object and can be computed as
\begin{verbatim}
  f = ballfun(@(x,y,z) cos(x.*y));% ballfun representing cos(xy)
  diff(f, 1)                      % Compute ballfun representing df/dx
  ans = 
   ballfun object
       domain          r    lambda    theta
       unit ball         24        43       40
\end{verbatim}
Ballfun calls the constructor after each operation to readjust the grid sizes (see \cref{sec_determination_size}). Here, a discretization size of $24\times 43\times 40$ was determined necessary to resolve $\partial f/\partial x$ while $f$ is represented by a $21\times 41\times 37$ CFF series.

\subsection{Fast rotation algorithm using a nonuniform Fourier transform}\label{sec:rotate}

Rotating functions defined on the ball has applications in many fields, including quantum mechanics, inverse scattering, and geophysics. Ballfun has a \texttt{rotate} command to efficiently perform rigid-body rotations of functions. Every rigid-body rotation can be specified by an Euler angle $(\alpha, \beta, \gamma)$ in the Z-X-Z convention~\cite{arfken1985mathematical}, which corresponds to rotating first by $\alpha$ around the $z$-axis then rotating by $\beta$ around the (original) $x$-axis and then, finally, rotating by $\gamma$ around the new $z$-axis. All the angles are given in radians. The algorithm to achieve this rotation requires a nontrivial computation because the rotated function must be represented by an approximant in the original coordinate system. 

The classical algorithm for computing the rotation of a function $f$ on the ball is to first express $f$ in terms of a spherical harmonic expansion and then to use the fact that the spherical harmonics form a basis of $SO(3)$~\cite{gelfand2018representations}.  Since Ballfun does not represent functions using spherical harmonic expansions, we use an algorithm based on the DFS method and the 2D nonuniform FFT~\cite{ruiz2018nonuniform}. We do this by taking the CFF grid in~\eqref{eq:grid} and rotating it by Euler angle $(\alpha,\beta,\gamma)$. Then, we evaluate the function at this rotated grid and call the Ballfun constructor. Since the rotated grid is almost always non-uniform in the $\theta$ and $\lambda$ variables of the doubled-up spherical coordinates and a Chebyshev grid in $r$, the evaluation is done in $\mathcal{O}(n^3\log n)$ operations with a 2D nonuniform FFT in $\theta$ and $\lambda$ and a DCT in $r$. 

\begin{figure}[htbp]
\centering
\begin{minipage}{.49\textwidth} 
\centering
\begin{overpic}[width=\textwidth,trim={120 80 110 80},clip]{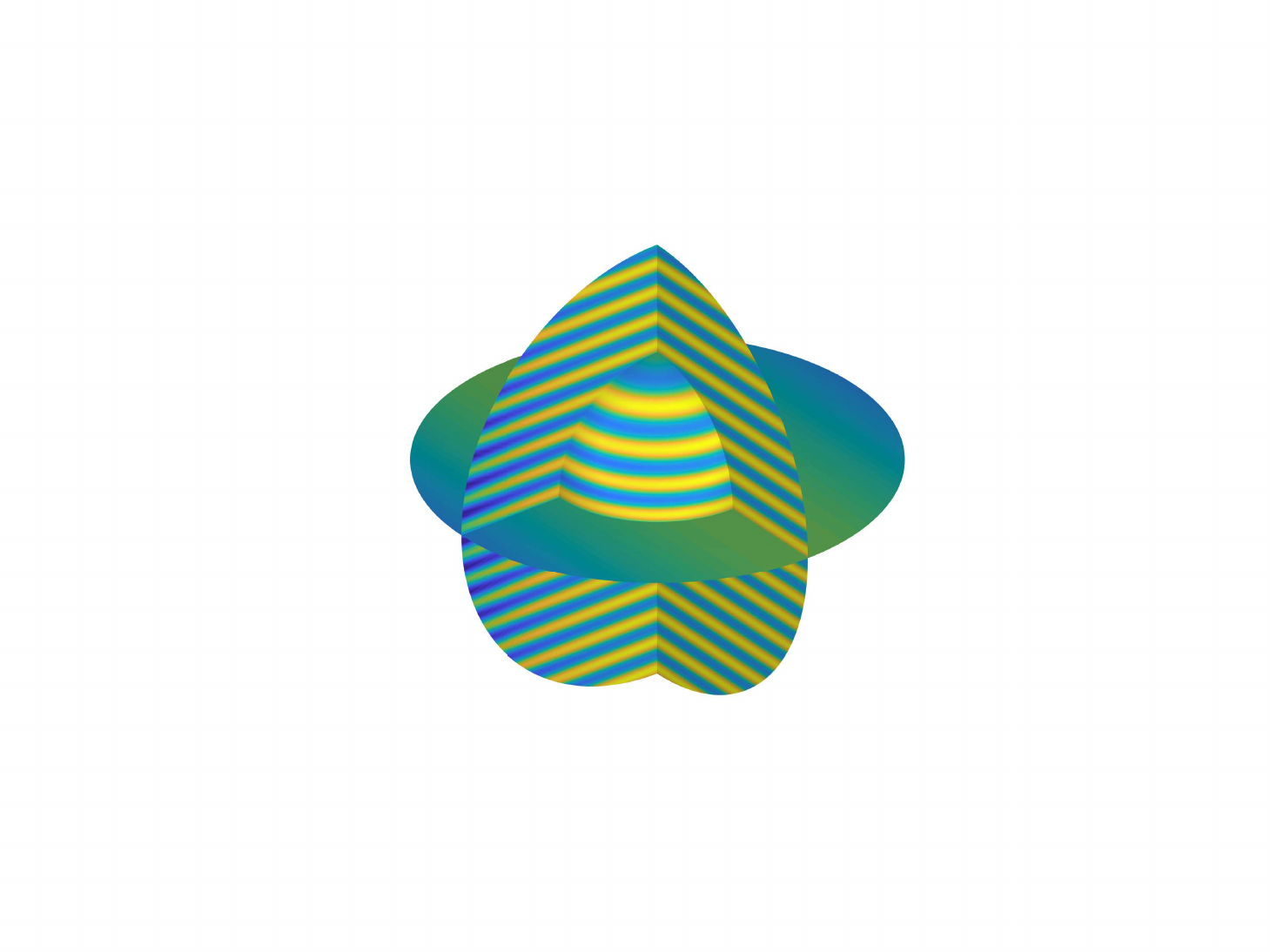}
\end{overpic} 
\end{minipage}
\begin{minipage}{.49\textwidth} 
\centering
\begin{overpic}[width=\textwidth,trim={120 80 110 80},clip]{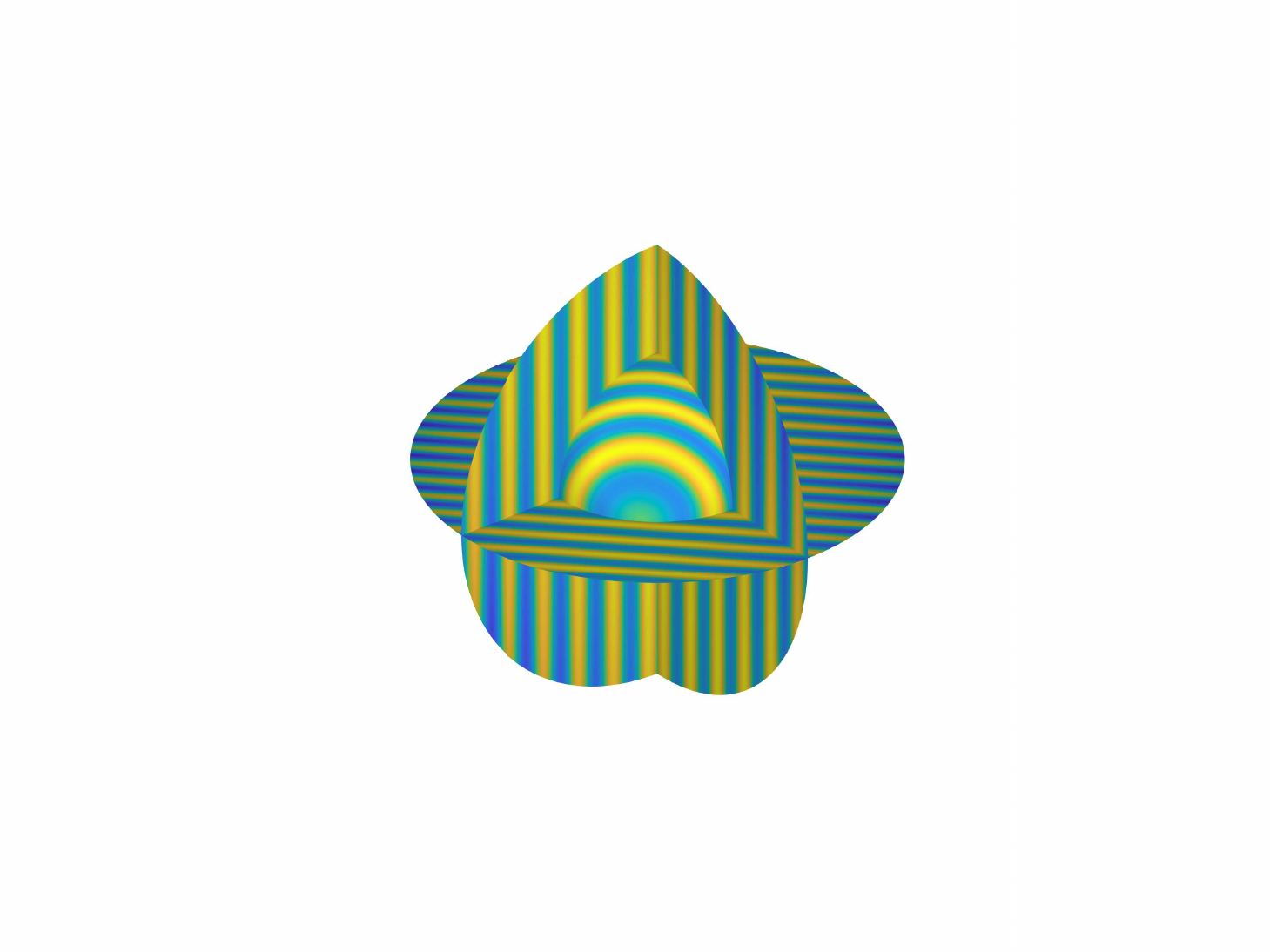}
\end{overpic} 
\end{minipage}
\label{fig_Rotate}
\caption{The function $f(x,y,z) = \sin(50z) - x^2$ (left) and its rotation by the Euler angles $\phi = -\pi/4$, $\theta = \pi/2$, and $\psi = \pi/8$ (right). The rotation of the function is computed with the \texttt{rotate} command and the functions are visualized with the \texttt{plot} command.}
\end{figure}

For example, the following code snippet calculates the rotation of $f(x,y,z) = \sin(50z) - x^2$ by Euler angle $(-\pi/4,\pi/2,\pi/8)$ (see~\cref{fig_Rotate}). 
\begin{verbatim} 
f = ballfun(@(x,y,z) sin(50*z) - x.^2) % ballfun for sin(50z)-x^2
f = 
   ballfun object
       domain          r    lambda    theta
       unit ball       90         5        179
g = rotate(f, -pi/4, pi/2, pi/8)       % Rotate by (-pi/4,pi/2,pi/8)
g = 
   ballfun object
       domain          r    lambda    theta
       unit ball       91      180       182
\end{verbatim} 
As one can see, the rotate command is also adaptive and selects the appropriate discretization to resolve the rotated function. 

\section{Fast spectral method for solving Helmholtz equation}\label{sec_Helmholtz}
In this section, we describe a fast algorithm for solving the Helmholtz equation on the ball with Neumann boundary conditions. An optimal-complexity algorithm for solving Helmholtz equation with Dirichlet conditions on the ball is described in~\cite{fortunato2017fast}, though it cannot immediately be generalized to the situation with Neumann conditions. Helmholtz solvers are useful in computational fluid dynamics as well as the computation of vector decompositions such as the poloidal-toroidal and Helmholtz--Hodge decompositions~\cite{backus1986poloidal,bhatia2013helmholtz}. 

\subsection{Discretization of Helmholtz equation}\label{sec_discretization_helmholtz}
Consider the Helmholtz equation on the ball, i.e., $u_{xx}+u_{yy}+u_{zz}+K^2u=f$ with Neumann boundary conditions $g(x,y,z)$ on the sphere $x^2+y^2+z^2=1$ and a real wave number $K$. The change of variables given by $(x,y,z)=(r\cos\lambda\sin\theta,r\sin\lambda\sin\theta,r\cos\theta)$, where $r\in [0,1]$, $\lambda\in[-\pi,\pi]$, and $\theta\in[0,\pi]$, transforms the equation into
\begin{equation}
\label{eq_helmholtz}
\frac{1}{r^2}\frac{\partial}{\partial r}\left(r^2\frac{\partial u}{\partial r}\right) + \frac{1}{r^2\sin\theta}\frac{\partial}{\partial\theta}\left(\sin\theta\frac{\partial u}{\partial\theta}\right)+\frac{1}{r^2\sin^2\theta}\frac{\partial^2u}{\partial\lambda^2} + K^2 u = f.
\end{equation}
One can multiply~\cref{eq_helmholtz} by $r^2\sin^2\theta$ to remove the singularities in the variable coefficients at the origin and at the poles of the ball. We then use the DFS method (see~\cref{sec_DFS}) to represent $u$ over the domain $(r,\lambda,\theta)\in[-1,1]\times[-\pi,\pi]\times[-\pi,\pi]$. This allows us to solve~\cref{eq_helmholtz} by seeking a tensor of Chebyshev--Fourier--Fourier coefficients for $u$ (see~\cref{eq:expansion}). 

Let $U=(u_{ijk})$ and $F=(f_{ijk})$ be $m\times n\times p$ tensors of CFF coefficients of $u$ and $(r^2\sin^2\theta) f$, respectively. Since~\cref{eq_helmholtz} decouples in the azimuthal variable $\lambda$, the following equation holds for $-n/2\leq j\leq n/2-1$:
\begin{equation}
\label{eq_helmholtz_2}
\sin^2\theta\frac{\partial}{\partial r}\left(r^2\frac{\partial u_j}{\partial r}\right) + \sin\theta\frac{\partial}{\partial\theta}\left(\sin\theta\frac{\partial u_j}{\partial\theta}\right) + \left(K^2 r^2\sin^2\theta-j^2\right)u_j = f_j,
\end{equation}
where the functions $u_j$ and $f_j$ are defined by
\[u_j(r,\theta)=\sum_{i=0}^{m-1}\sum_{k=-p/2}^{p/2-1} u_{ijk} T_{i}(r)e^{\mathbf{i}k\theta},\quad\quad f_j(r,\theta)=\sum_{i=0}^{m-1}\sum_{k=-p/2}^{p/2-1} f_{ijk} T_{i}(r)e^{\mathbf{i}k\theta}.\]
We discretize~\cref{eq_helmholtz_2} in the radial variable using the ultraspherical spectral method~\cite{olver2013fast}, and in the polar variable using the Fourier spectral method. Partial derivatives in the polar variable $\theta$ and multiplication by $\sin\theta$ are represented by sparse and banded matrices in the Fourier basis. The ultraspherical spectral method results in sparse and banded matrices of operators (such as differentiation or multiplication by $r$) between Chebyshev and ultraspherical polynomials. This allows us to write~\cref{eq_helmholtz_2} in the form of a generalized Sylvester equation~\cite{gardiner1992solution} in the unknown matrix $U(:,j,:)$:
\begin{equation}
\label{eq_Helmholtz_Sylvester}
L_rU(:,j,:)M_{\sin^2}^\top+S_{02}U(:,j,:)L_\theta^\top=S_{02}F(:,j,:),
\end{equation}
where $L_r$ is the matrix representing the operator $u\mapsto\frac{\partial u}{\partial r}\left(r^2\frac{\partial u}{\partial r}\right)+K^2r^2u$ from the Chebyshev basis $T$ to the ultraspherical basis $C^{(2)}$ and $S_{02}$ is the conversion matrix between these bases~\cite{olver2013fast}. The matrices $M_{\sin^2}$ and $L_\theta$ represent the multiplication by $\sin^2\theta$ and the operator $u\mapsto\sin\theta\frac{\partial u}{\partial\theta}\left(\sin\theta\frac{\partial u}{\partial\theta}\right)-j^2u$ in the Fourier basis, respectively.

\subsection{Imposing Neumann boundary conditions when $\mathbf{K\neq 0}$}
It is essential to slightly modify~\cref{eq_Helmholtz_Sylvester} to impose Neumann boundary conditions on $u$, i.e., $\partial_ru|_{r=1}=g(\lambda,\theta)$.  The first step is to double-up the smooth function  $g(\lambda,\theta)$ in the $\theta$ variable using the DFS method~\cite{townsend2016computing} and define its Fourier--Fourier matrix of coefficients $G^+=(g^+_{jk})$. Since the radial variable $r$ of~\cref{eq_helmholtz} has been doubled-up, we need to impose a Neumann condition at $r=1$ and $r=-1$. The matrix of coefficients of the boundary condition at $r=-1$, $G^-=(g^-_{jk})$, can be deduced from $G^+$ (see \cref{sec_DFS}) and takes the form
\[G^-(j,:) = (-1)^jG^+(j,:), \quad -\frac{n}{2}\leq j\leq \frac{n}{2}-1.\]
The Neumann operators $u\mapsto\left.\frac{\partial u}{\partial r}\right|_{r= 1}$ and $u\mapsto\left.\frac{\partial u}{\partial r}\right|_{r= -1}$ are represented in the Chebyshev basis by the $1\times m$ matrices 
\[B^+(k) = k^2,\quad B^-(k) = (-1)^{k+1}k^2,\quad 0\leq k\leq m-1,\]
respectively.
The Neumann conditions also decouple in the variable $\lambda$ and can be written as
\begin{equation}
\label{eq_Neumann}
\begin{pmatrix}B^+\\B^-\end{pmatrix}U(:,j,:)=
\begin{pmatrix}G^+(j,:)\\G^-(j,:)\end{pmatrix},\quad -\frac{n}{2}\leq j\leq \frac{n}{2}.
\end{equation}
Therefore, a Helmholtz's solver~\cref{eq_helmholtz} with Neumann boundary conditions is realized by solving the following $m\times p$ generalized Sylvester equation with linear constraints:
\begin{align}
L_rU(:,j,:)M_{\sin^2}^\top+S_{02}U(:,j,:)L_\theta^\top=S_{02}F(:,j,:),\label{eq_Sylvester_1}
\\
\begin{pmatrix}B^+\\B^-\end{pmatrix}U(:,j,:)=
\begin{pmatrix}G^+(j,:)\\G^-(j,:)\end{pmatrix},\label{eq_Sylvester_2}
\end{align}
where $-n/2\leq j\leq n/2-1$. The constraints~\cref{eq_Sylvester_2} can be used to remove degrees of freedom in $U(:,j,:)$ and transform~\cref{eq_Sylvester_1} into a generalized Sylvester equation with a unique solution without constraints~\cite{townsend2015automatic}, i.e.,
\begin{equation}
\label{eq_modified_sylvester}
\tilde{L}_rX_jM_{\sin^2}^\top+\tilde{S}_{02}X_jL_\theta^\top=\tilde{S}_{02}\tilde{F}(:,j,:), \quad -\frac{n}{2}\leq j\leq \frac{n}{2}-1.
\end{equation}

\begin{figure}[htbp]
\begin{minipage}{.49\textwidth} 
\centering
\begin{overpic}[width=\textwidth,trim={40 0 40 0},clip]{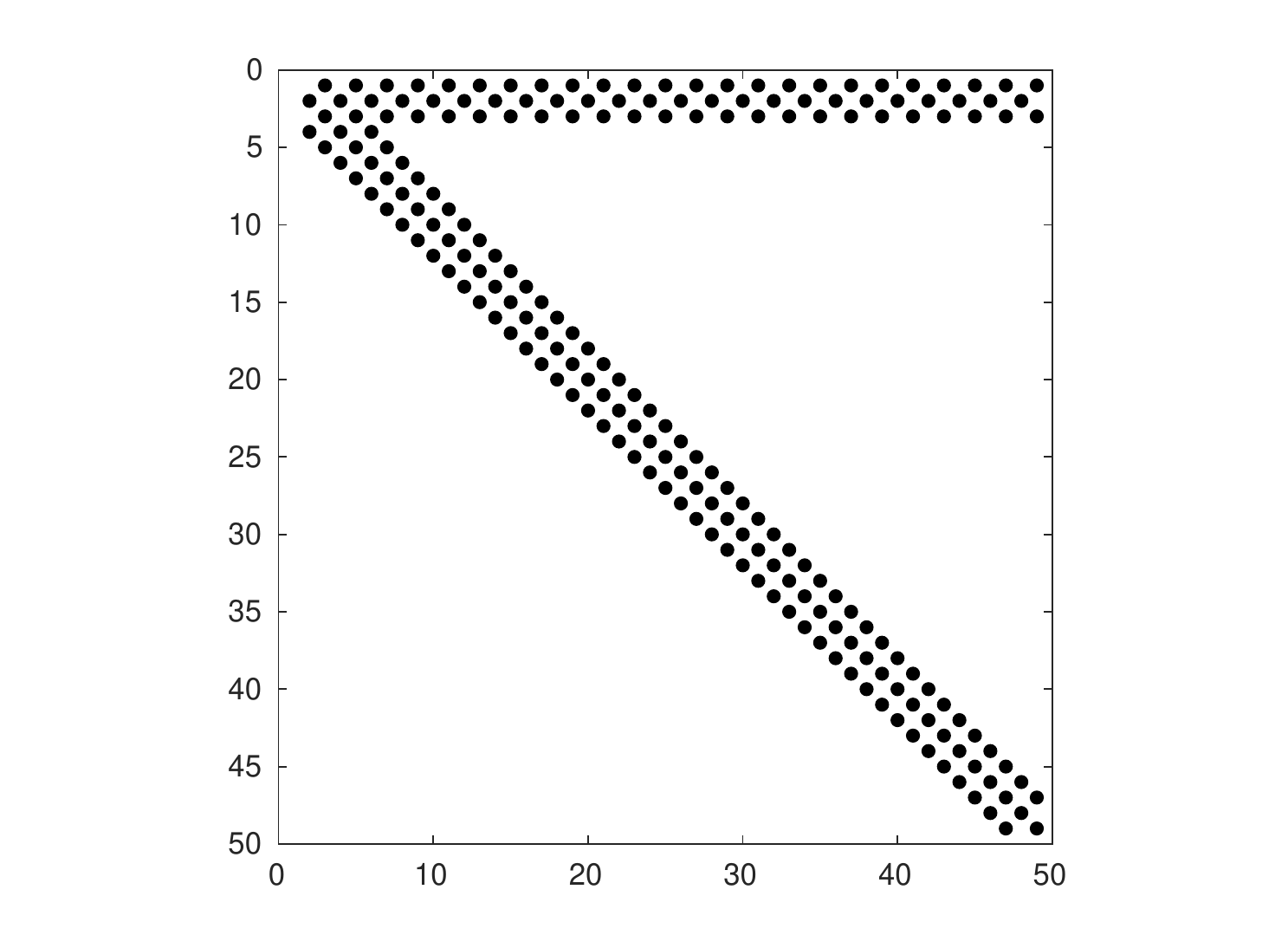}
\end{overpic} 
\end{minipage}
\begin{minipage}{.49\textwidth} 
\centering
\begin{overpic}[width=\textwidth,trim={40 0 40 0},clip]{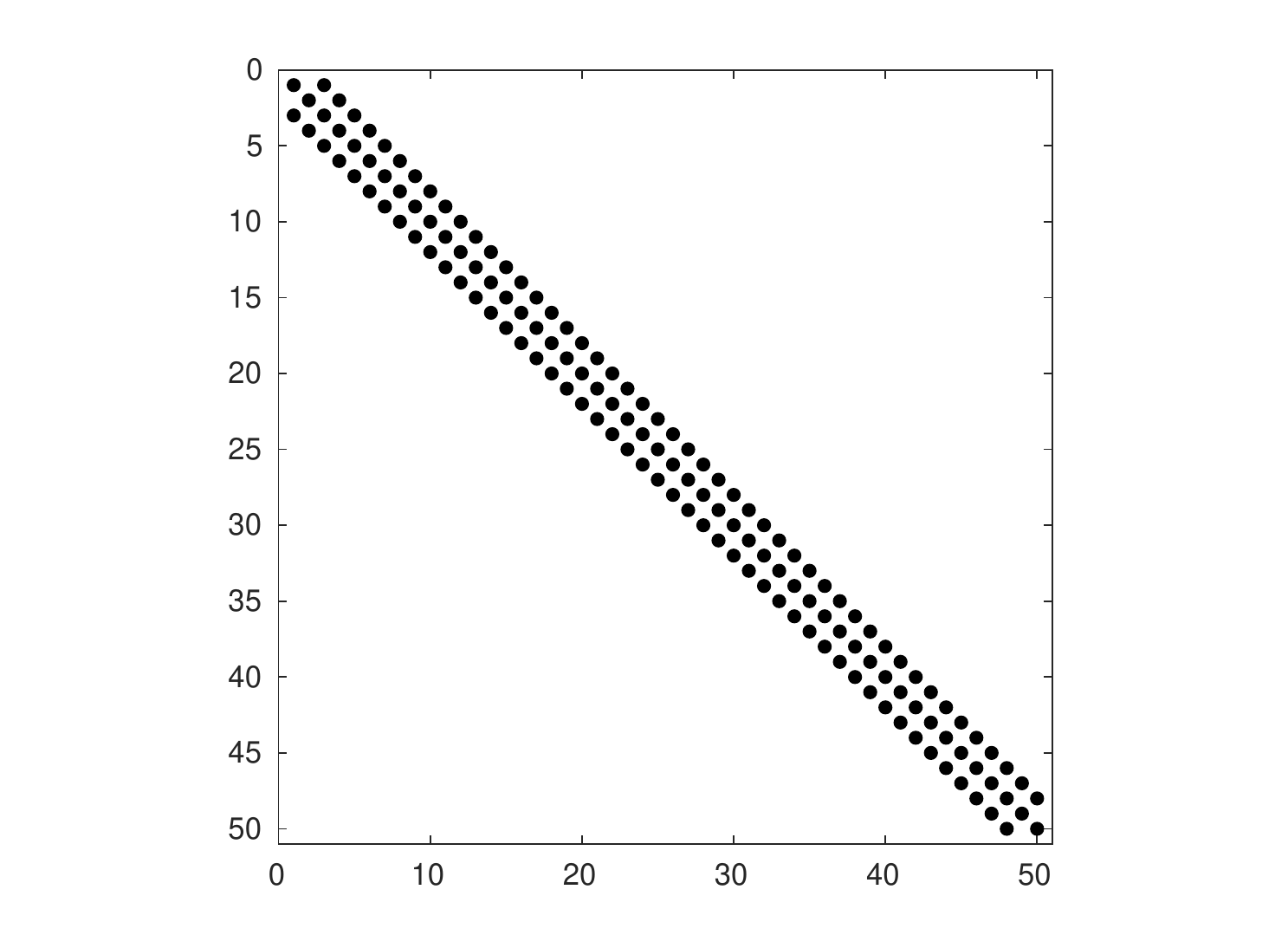}
\end{overpic} 
\end{minipage}
\label{fig_matrix_Neumann}
\caption{Sparsity structure of the matrix $\tilde{L}_r$ (left) and sparsity structure of the banded matrix $M_{\sin^2}$ (right) for $m=p=50$.}
\end{figure}

\cref{fig_matrix_Neumann} shows the sparsity structure of the matrices $\tilde{L}_r$ and $M_{\sin^2}$ in~\cref{eq_modified_sylvester}. We obtain $n$ decoupled Sylvester matrix equations, where each one can be solved in $\mathcal{O}(m^3+p^3)$ operations using the Bartels--Stewart algorithm~\cite{bartels1972solution,gardiner1992solution}. Once $X_j$ has been computed, the matrix of coefficients $U(:,j,:)$ can be recovered using the linear constraints.  Thus, the total complexity is $\mathcal{O}((m^3+p^3)n)$ operations. 

\subsection{Imposing Neumann boundary conditions when $\mathbf{K = 0}$}\label{sec_0_Neumann}
We consider the zeroth Fourier mode $j=0$ of~\cref{eq_Helmholtz_Sylvester} with $K=0$ (Poisson equation). The solution to this equation with Neumann boundary conditions is unique only up to a constant. However, this additional constraint cannot be imposed on a Sylvester matrix equation. Therefore, we transform~\cref{eq_Helmholtz_Sylvester} into the Chebyshev--Legendre basis to decouple this Sylvester equation in the polar variable $\theta$. The function $u_0(r,\theta)$, defined in \cref{sec_discretization_helmholtz}, satisfies the following equation:
\begin{equation}
\label{eq_Poisson_0}
\frac{\partial}{\partial r}\left(r^2\frac{\partial u_0}{\partial r}\right) + \frac{1}{\sin\theta}\frac{\partial}{\partial\theta}\left(\sin\theta\frac{\partial u_0}{\partial\theta}\right)= r^2f_0.
\end{equation}
The functions $u_0$ and $r^2f_0$ can be expressed in the Chebyshev--Legendre (CL) basis as
\[u_0(r,\theta) = \sum_{i=0}^{m-1}\sum_{k=0}^{p-1}\tilde{u}_{i0k}T_i(r)P_k(\cos\theta),\quad r^2f_0(r,\theta) = \sum_{i=0}^{m-1}\sum_{k=0}^{p-1}\tilde{f}_{i0k}T_i(r)P_k(\cos\theta).\]

The zeroth Fourier modes $j=0$ of the Neumann boundary conditions at $r=1$, $g^+_0(\theta)$, and at $r=-1$,  $g^-_0(\theta)$, can also be written as Legendre series
\[g^+_0(\theta) = \sum_{k=0}^{p-1}\tilde{g}^+_{0k}P_k(\cos\theta), \quad g^-_0(\theta) = \sum_{k=0}^{p-1}\tilde{g}^-_{0k}P_k(\cos\theta).\]
The orthogonality of the Legendre basis allows us to decouple~\cref{eq_Poisson_0} in the polar variable $\theta$ as $p$ ordinary differential equations with Neumann boundary conditions:
\begin{align}
\frac{\partial}{\partial r}\left(r^2\frac{\partial \tilde{u}_{0k}}{\partial r}\right) -k(k+1)\tilde{u}_{0k} &= \tilde{f}_{0k},\label{eq_Neumann_Legendre_1} \\
\left.\frac{\partial\tilde{u}_{0k}}{\partial r}\right|_{r=1}=\tilde{g}^+_{0k},\quad
\left.\frac{\partial\tilde{u}_{0k}}{\partial r}\right|_{r=-1}&=\tilde{g}^-_{0k},\label{eq_Neumann_Legendre_2}
\end{align}
for $0\leq k\leq p-1$. The functions $\tilde{u}_{0k}$ and $\tilde{f}_{0k}$ are defined by
\[\tilde{u}_{0k}(r) = \sum_{i=0}^{m-1}\tilde{u}_{i0k}T_i(r),\qquad \tilde{f}_{0k}(r) = \sum_{i=0}^{m-1}\tilde{f}_{i0k}T_i(r),\quad r\in [-1,1].\]
For each $0< k\leq p-1$,~\cref{eq_Neumann_Legendre_1,eq_Neumann_Legendre_2} can be solved in $\mathcal{O}(m)$ operations using the ultraspherical spectral method~\cite{olver2013fast}.
The case $k = 0$ is solved by the same technique with the additional linear constraint $\tilde{u}_{000} = 0$ to impose uniqueness of the global solution $u$.

Once the matrix of Chebyshev--Legendre coefficients $\tilde{U}_0=(\tilde{u}_{i0k})$ has been computed, it can be converted to the Chebyshev--Fourier basis in $\mathcal{O}(mp\log p)$ operations using the Legendre--Chebyshev tranform~\cite{potts1998fast}. 

\subsection{Numerical examples} \label{sec_ex_helmholtz}

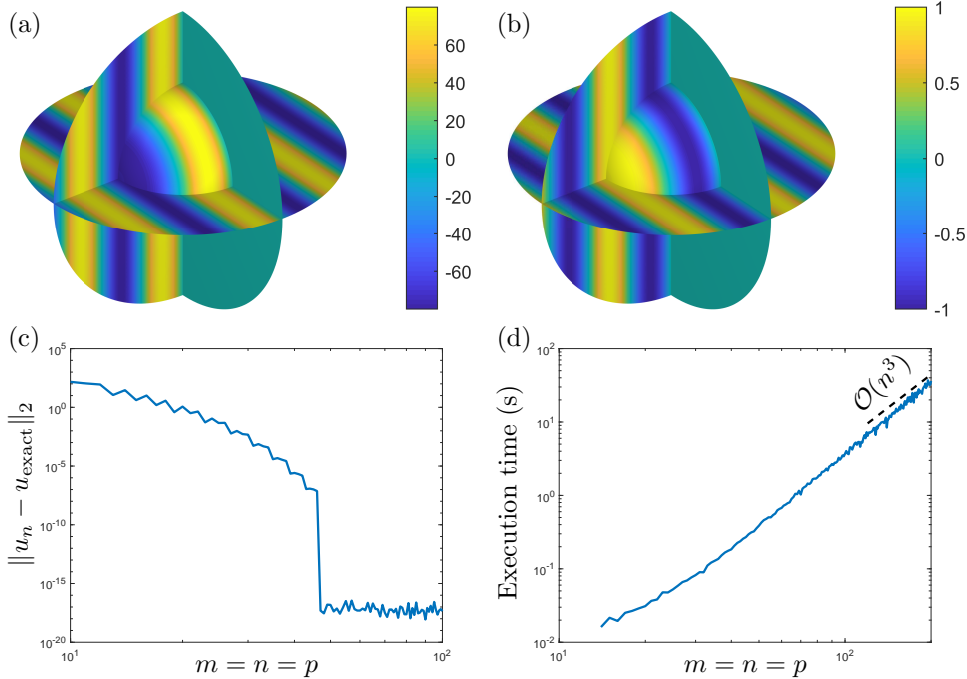
\begin{figure}[htbp]
\begin{minipage}{.49\textwidth} 
\begin{overpic}[width=\textwidth,trim={130 80 50 75},clip]{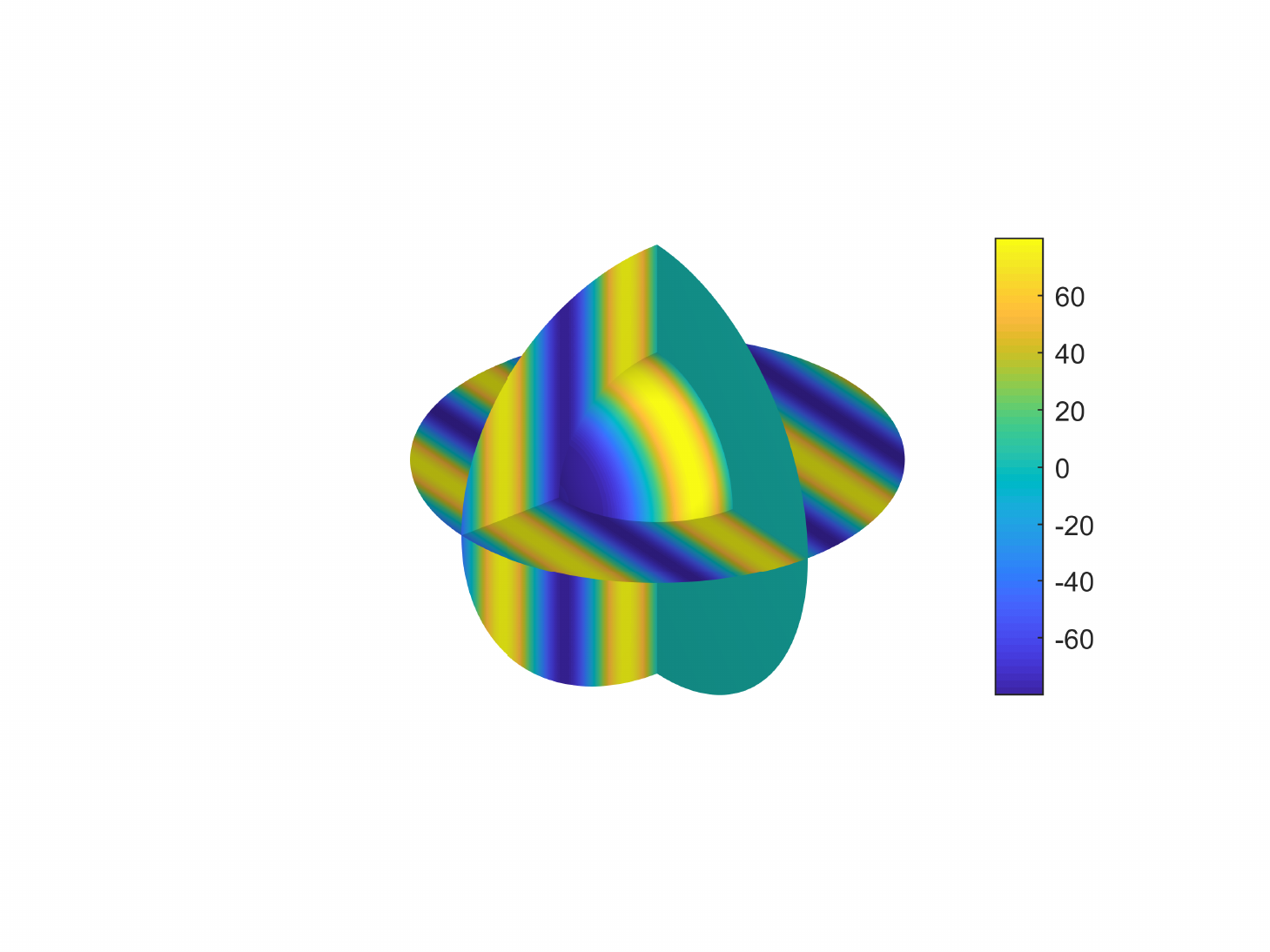}
\put(0,60){(a)}
\end{overpic} 
\end{minipage}
\begin{minipage}{.49\textwidth} 
\begin{overpic}[width=\textwidth,trim={130 80 50 75},clip]{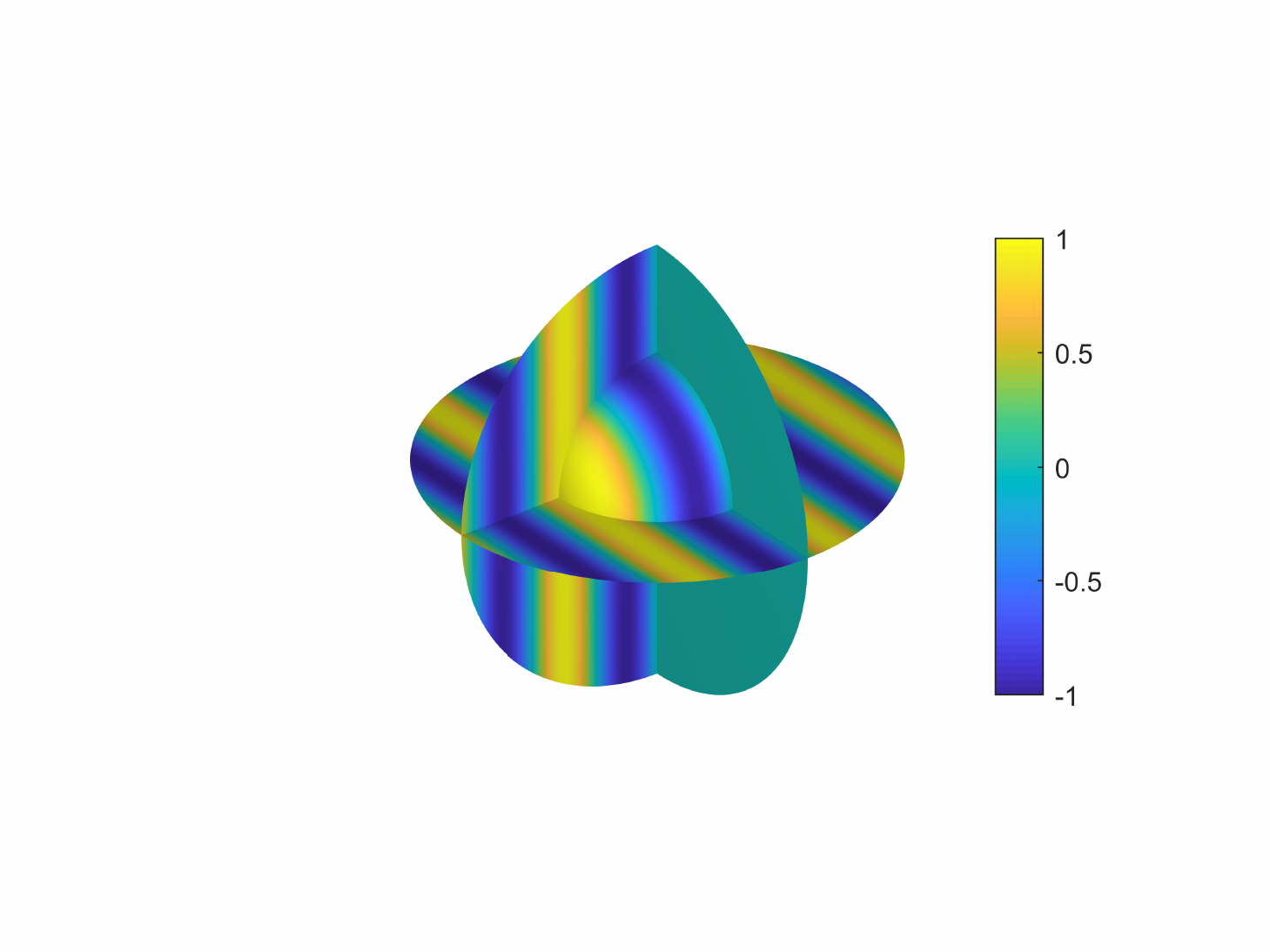}
\put(0,60){(b)}
\end{overpic} 
\end{minipage}\\
\begin{minipage}{.49\textwidth} 
\begin{overpic}[width=\textwidth]{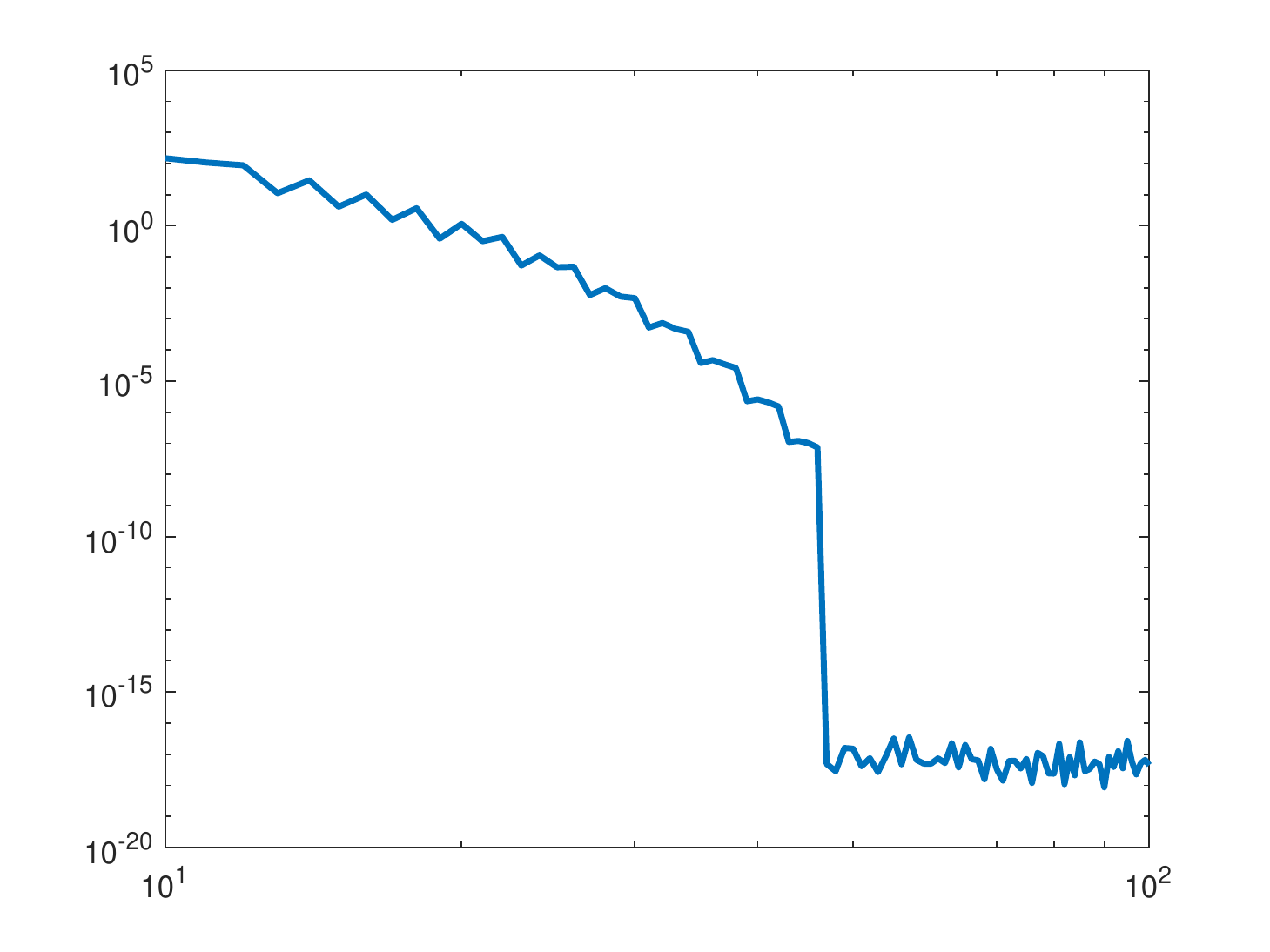}
\put(0,70){(c)}
\put(39,2){$m=n=p$}
\put(0,24){\rotatebox{90}{$\|u_n-u_{\text{exact}}\|_2$}}
\end{overpic} 
\end{minipage}
\begin{minipage}{.49\textwidth} 
\begin{overpic}[width=\textwidth]{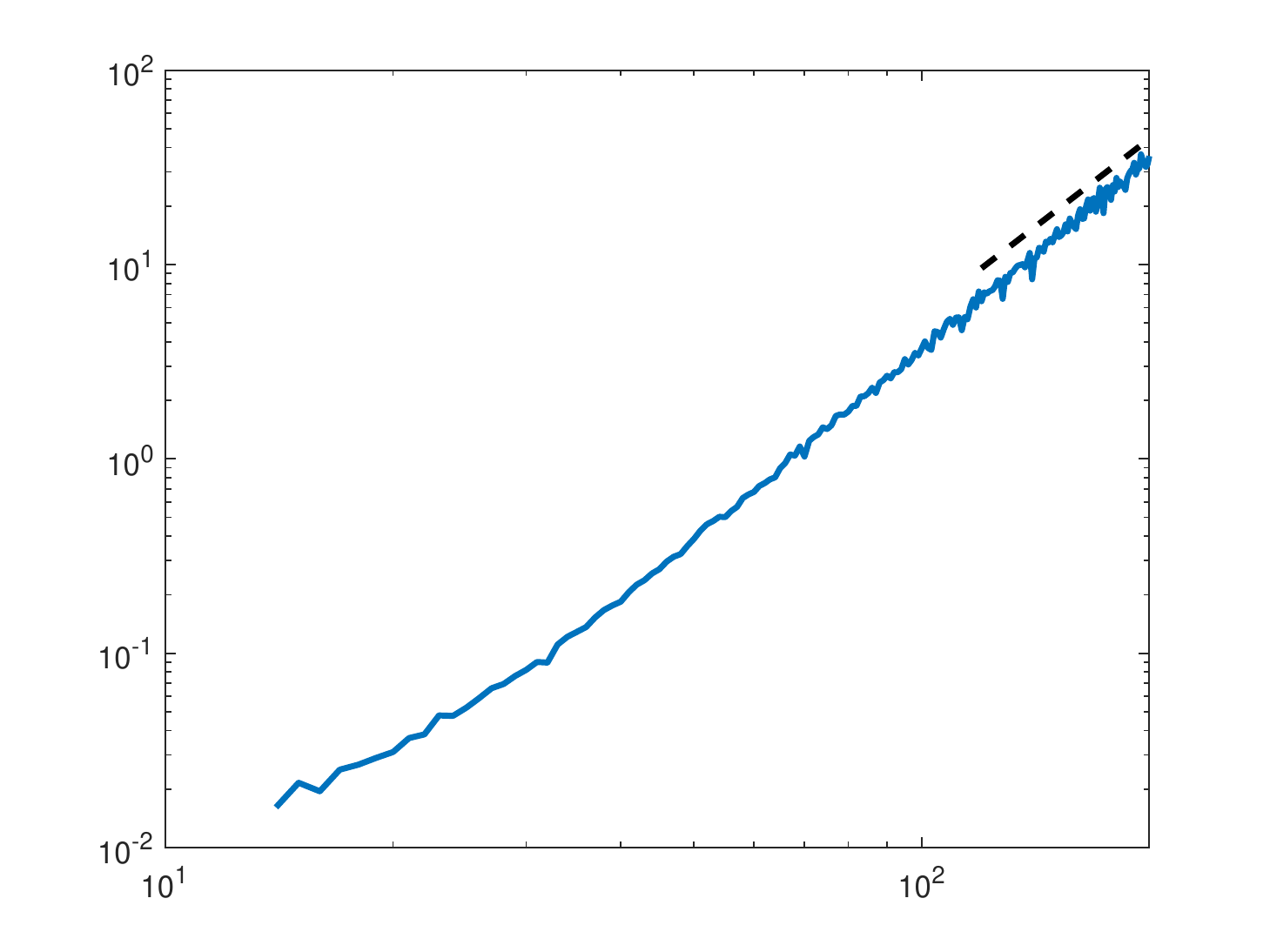}
\put(0,70){(d)}
\put(39,2){$m=n=p$}
\put(0,17){\rotatebox{90}{Execution time (s)}}
\put(73,56){\rotatebox{45}{$\mathcal{O}(n^3)$}}
\end{overpic} 
\end{minipage}
\label{fig_Helmholtz_Neumann}
\caption{Function $f(x,y,z)=-80\sin(10x)$ (see (a)) and solution $u$ to the Helmholtz equation $\Delta u + 20u = f(x,y,z)$ with Neumann boundary conditions $g(x,y,z) = 10x\cos(10x)$ (see (b)). The solution $u$ is computed with a spectral discretization of $m=n=p=50$. Error in the 2-norm between the exact solution $u_\text{exact}=\sin(10x)$ and the computed solution obtained with the \texttt{helmholtz} command (see (c)). The computational timings (see (d)).}
\end{figure}

In \cref{fig_Helmholtz_Neumann} (a) we plot a solution to the Helmholtz equation
\begin{equation} \label{eq_helmholtz_numerics}
\nabla^2 u + 20u = -80\sin(10x),
\end{equation}
with Neumann boundary conditions $g(x,y,z)=10x\cos(10x)$. The error between the computed and the exact solution to \cref{eq_helmholtz_numerics} is shown in \cref{fig_Helmholtz_Neumann} (c) and confirms the spectral convergence of our method. The computed solution is resolved to machine precision for $n\geq 50$. Our Helmholtz solver on the ball can be invoked in Ballfun via the \texttt{helmholtz} command.
\begin{verbatim} 
  rhs = ballfun(@(x,y,z)-80*sin(10*x));    % Right-hand side
  bc = @(x,y,z)10*x.*cos(10*x);            % Boundary conditions
  K = sqrt(20);                            % Wave number
  n = 50;                                  % Spectral discretization
  u = helmholtz(rhs, K, bc, n, 'neumann'); % Helmholtz solver
\end{verbatim} 

The execution times\footnote{Timings were performed on a 3.3 GHz Intel Core i5 CPU using MATLAB 2018a without explicit parallelization.} to solve \cref{eq_helmholtz_numerics} for different discretization sizes $n$ (see~\cref{fig_Helmholtz_Neumann} (d)). We can then solve Helmholtz's equation on a ball with one million degrees of freedom in a few seconds on a standard CPU.

\begin{figure}[htbp]
\centering
\begin{minipage}{.49\textwidth} 
\begin{overpic}[width=\textwidth,trim={130 80 50 60},clip]{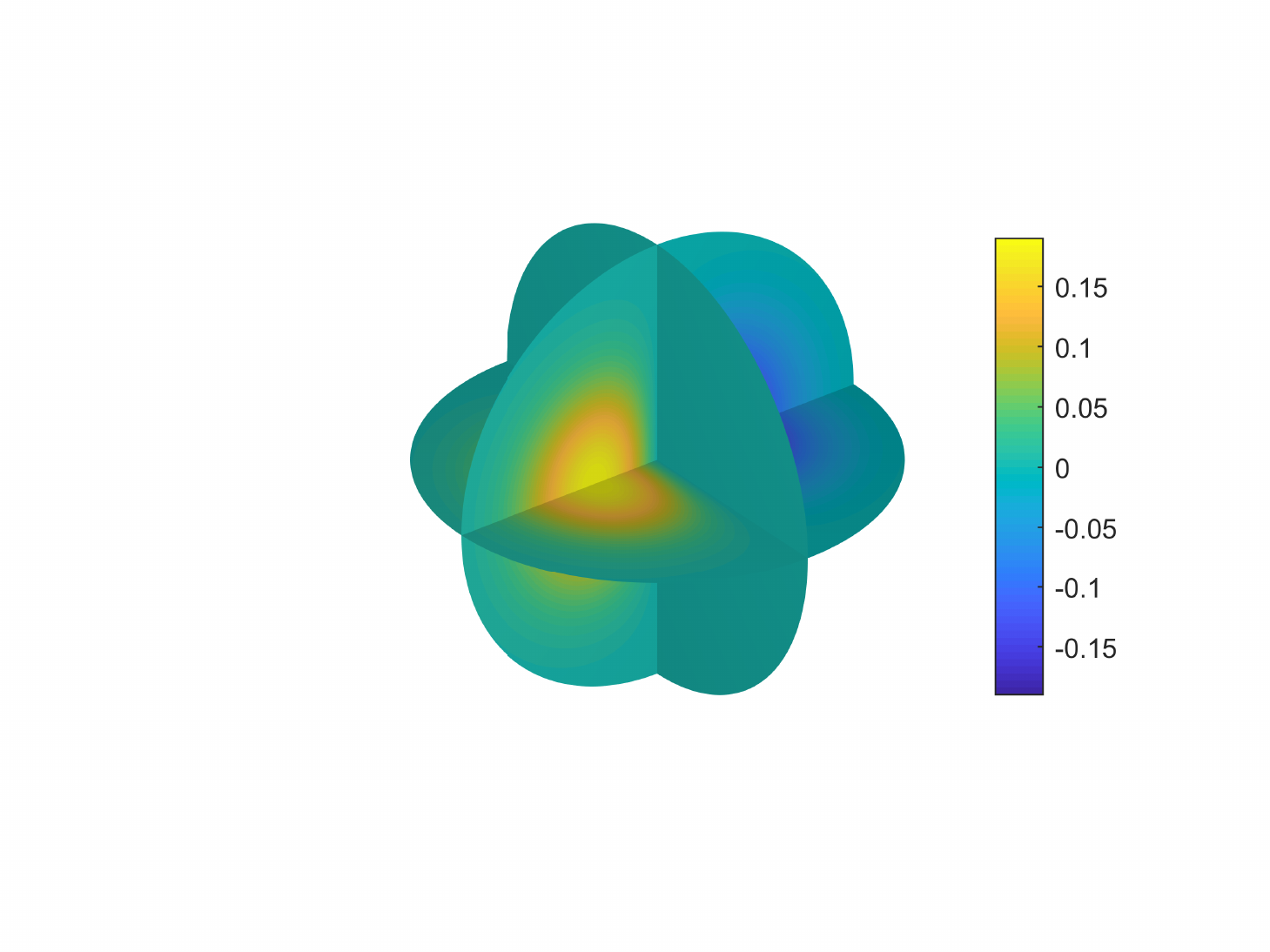}
\put(0,60){(a)}
\end{overpic}
\end{minipage}
\begin{minipage}{.49\textwidth} 
\begin{overpic}[width=\textwidth,trim={130 80 50 60},clip]{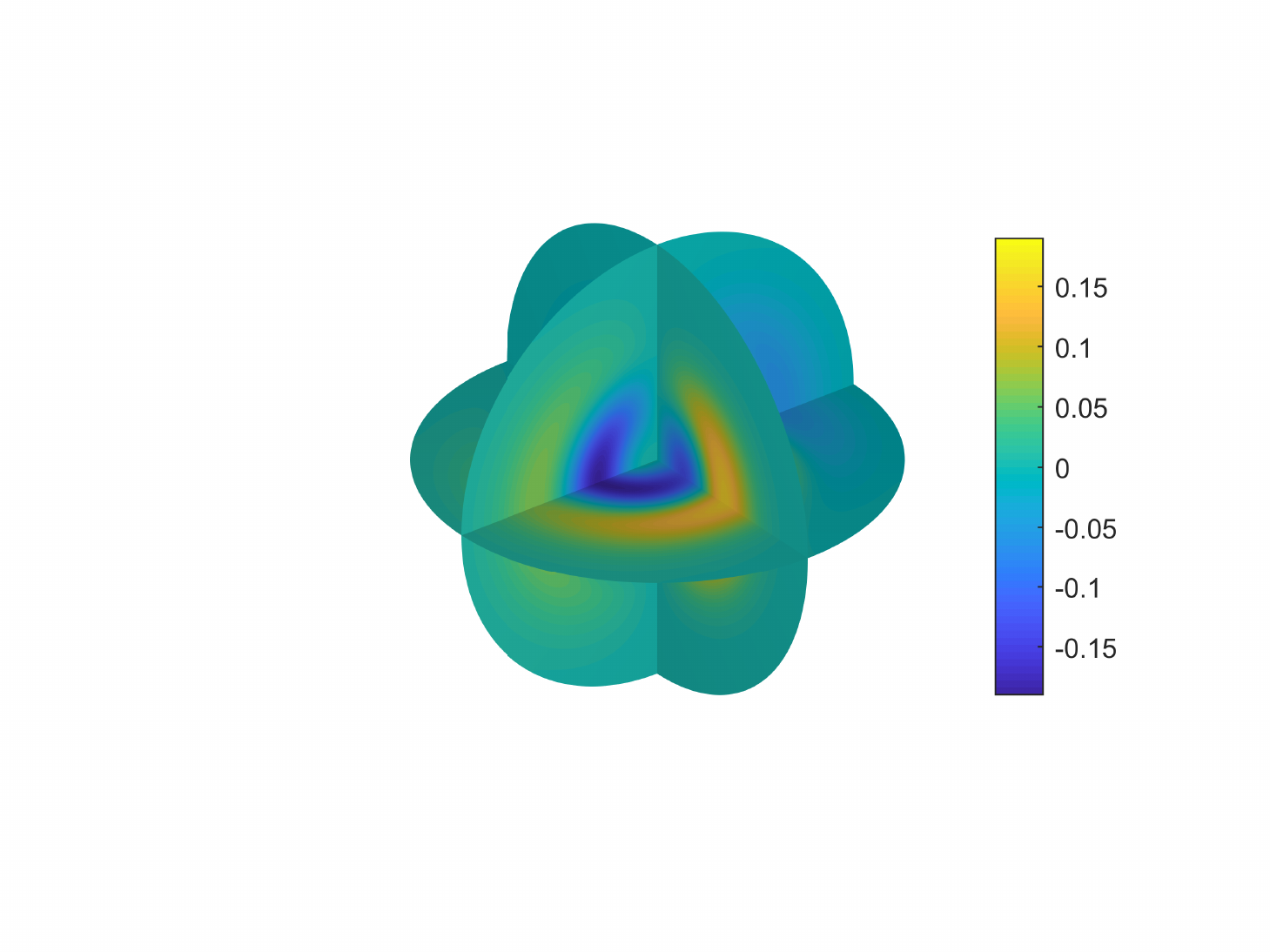}
\put(0,60){(b)}
\end{overpic} 
\end{minipage}\\
\begin{minipage}{.49\textwidth} 
\begin{overpic}[width=\textwidth,trim={130 80 50 60},clip]{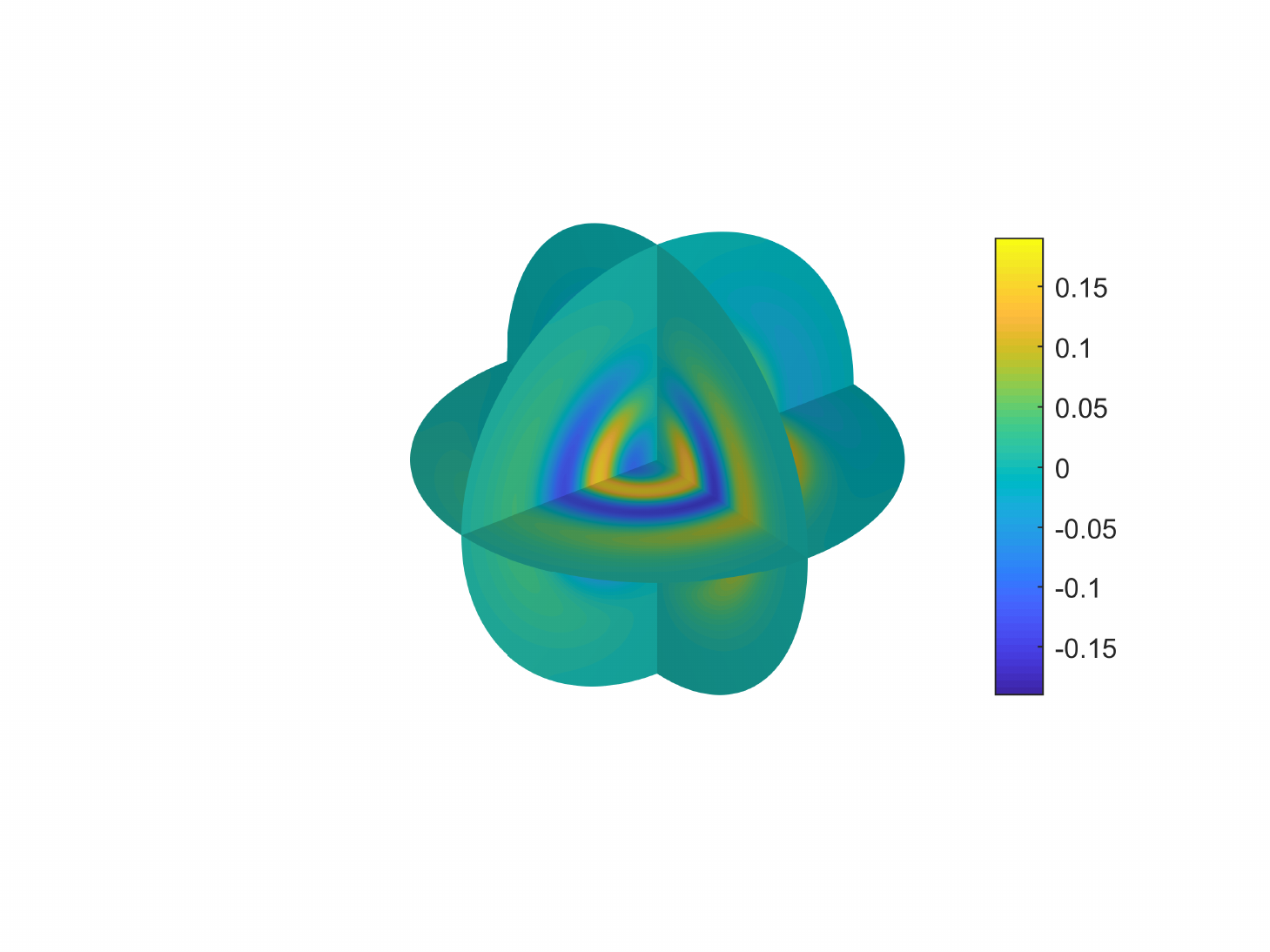}
\put(0,60){(c)}
\end{overpic} 
\end{minipage}
\begin{minipage}{.49\textwidth} 
\begin{overpic}[width=\textwidth,trim={130 80 50 60},clip]{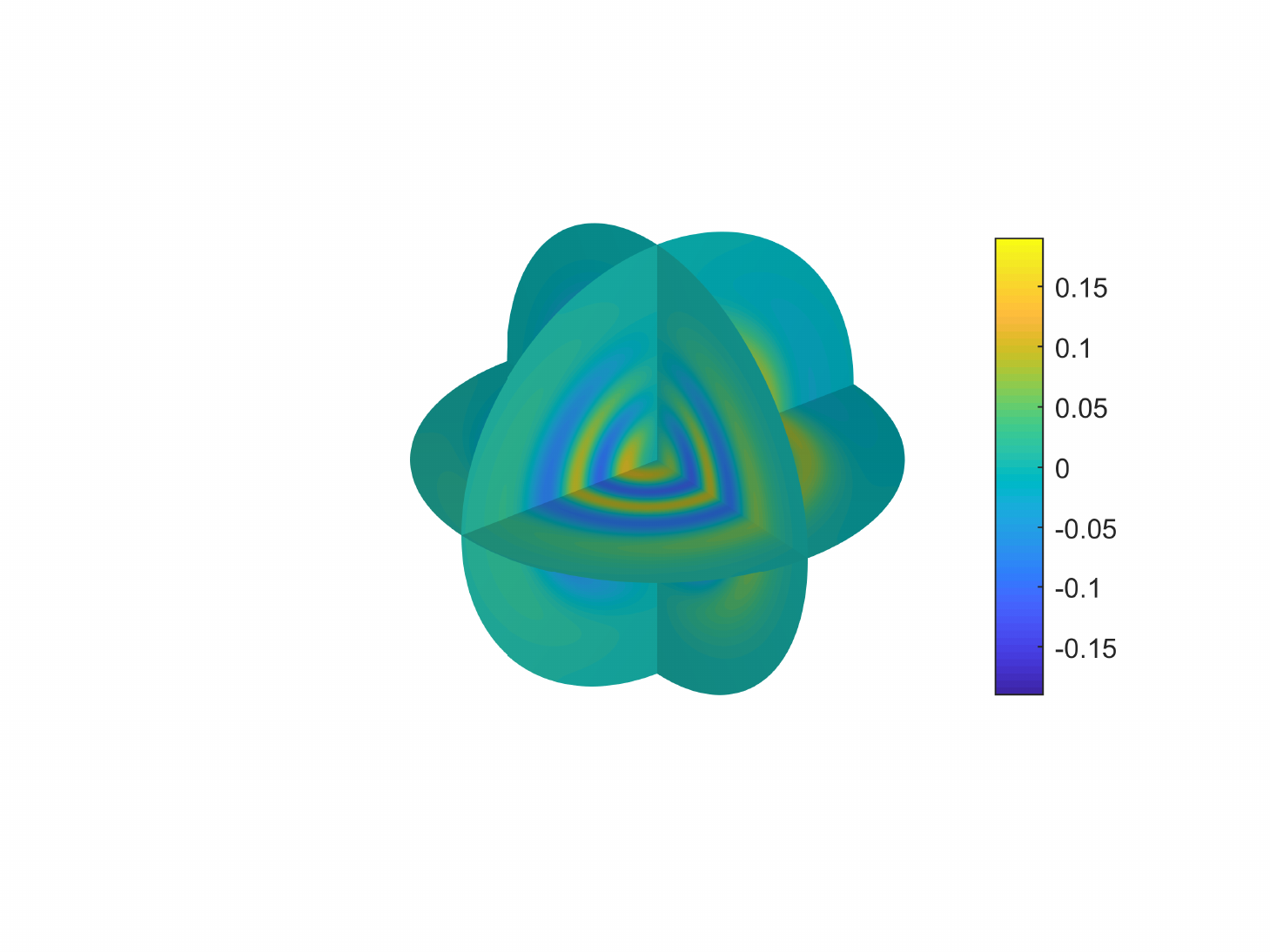}
\put(0,60){(d)}
\end{overpic} 
\end{minipage}
\label{fig_adv_diff}
\caption{Solutions to the advection-diffusion equation at $t = 0$ (see (a)), $t=5$ (see (b)), $t=10$ (see (c)), and $t=15$ (see (d)).}
\end{figure}

As a second example we consider the advection-diffusion equation on the unit ball
\begin{equation} \label{eq_adv_diff}
\frac{\partial c}{\partial t}=D\nabla^2c-v\cdot\nabla c,
\end{equation}
where $D$ is the diffusion coefficient and $v$ is a divergence-free vector field. We choose $D=1/5000$ and $v = \nabla\times[ze^{-5(x^2+y^2+z^2)}(x,y,z)]$ to satisfy the no-slip condition. The no-flux condition for $c$ reduces to $\partial c/\partial \vec{n}=0$ at the boundary. We impose the initial condition $c(x,y,z) = -xe^{-5(x^2+y^2+z^2)}$ and solve \cref{eq_adv_diff} by using the implicit-explicit backward differentiation of order one (IMEX-BDF1) scheme. This yields the following Helmholtz equation:
\[\nabla^2c^{n+1}+K^2c^{n+1}=K^2c^n+\frac{1}{D}v\cdot\nabla c^n,\quad \left.\frac{\partial c}{\partial \vec{n}}\right|_{\partial B(0,1)} = 0,\]
where $c_n$ denotes the solution at time $t = n\Delta t$, $\Delta t = 5\times 10^{-2}$ is the time step, and $K^2 = -1/(D\Delta t)$. The solution $c$ to \cref{eq_adv_diff} at different times is illustrated in \cref{fig_adv_diff}.

\section{Vector-valued functions on the ball}\label{sec_vector}
Ballfun is also designed to work with vector-valued functions defined on the unit ball as well as scalar-valued ones. Expressing vector-valued functions in spherical coordinates is inconvenient since the unit vectors in this coordinate system are singular at the poles of the ball~\cite{swarztrauber1981approximation}. Therefore, we express vector-valued functions in Cartesian coordinates as the components of the vector field are then themselves smooth functions. After using this convention, vector-valued functions introduce few complications from the point-of-view of approximation as each component is represented as an independent scalar function. A vector-valued function can be constructed in Ballfun as follows: 
\begin{verbatim}
V = ballfunv(@(x,y,z) sin(x), @(x,y,z) x.*y, @(x,y,z) cos(z))
ballfunv object containing
ballfun object:
    domain          r    lambda    theta
    unit ball      14      27       27
ballfun object:
    domain          r    lambda    theta
    unit ball       3       5        5
ballfun object:
    domain          r    lambda    theta
    unit ball      15       1       29
\end{verbatim}
It can be seen that a vector field is stored as a \texttt{ballfunv} object, which consists of three \texttt{ballfun} objects corresponding to the three components in Cartesian coordinates. Each \texttt{ballfun} object has its own discretization in $r$, $\lambda$ and $\theta$.

\subsection{Vector calculus on the ball}
The more interesting side of vector-valued functions in Ballfun is the set of operations that can be implemented, which are potentially useful for applications.   The standard operations for vector calculus such as the curl, the gradient or the divergence are implemented in Ballfun in the \texttt{curl}, \texttt{gradient}, and \texttt{divergence} commands, respectively. Due to the way we represent vector-valued functions, we compute these operations in Cartesian coordinate system. For example, the curl of a vector-valued function can be written as
\[
\nabla\times V = \left[ \frac{\partial V_z}{\partial y}-\frac{\partial V_y}{\partial z}, \frac{\partial V_x}{\partial z}-\frac{\partial V_z}{\partial x}, \frac{\partial V_y}{\partial x}-\frac{\partial V_x}{\partial y}\right]^T.
\]
The curl of $V(x,y,z)=(\sin x,xz,\cos z)$ can be computed and displayed using the \texttt{quiver} command (see~\cref{fig_Vector_calculus}):
\begin{verbatim}
  W = curl( V ); % Compute curl of V 
  quiver( W )    % Plot vector field W using quiver
\end{verbatim}

\begin{figure}[htbp]
\begin{minipage}{.49\textwidth} 
\begin{overpic}[width=\textwidth,trim={120 70 50 50},clip]{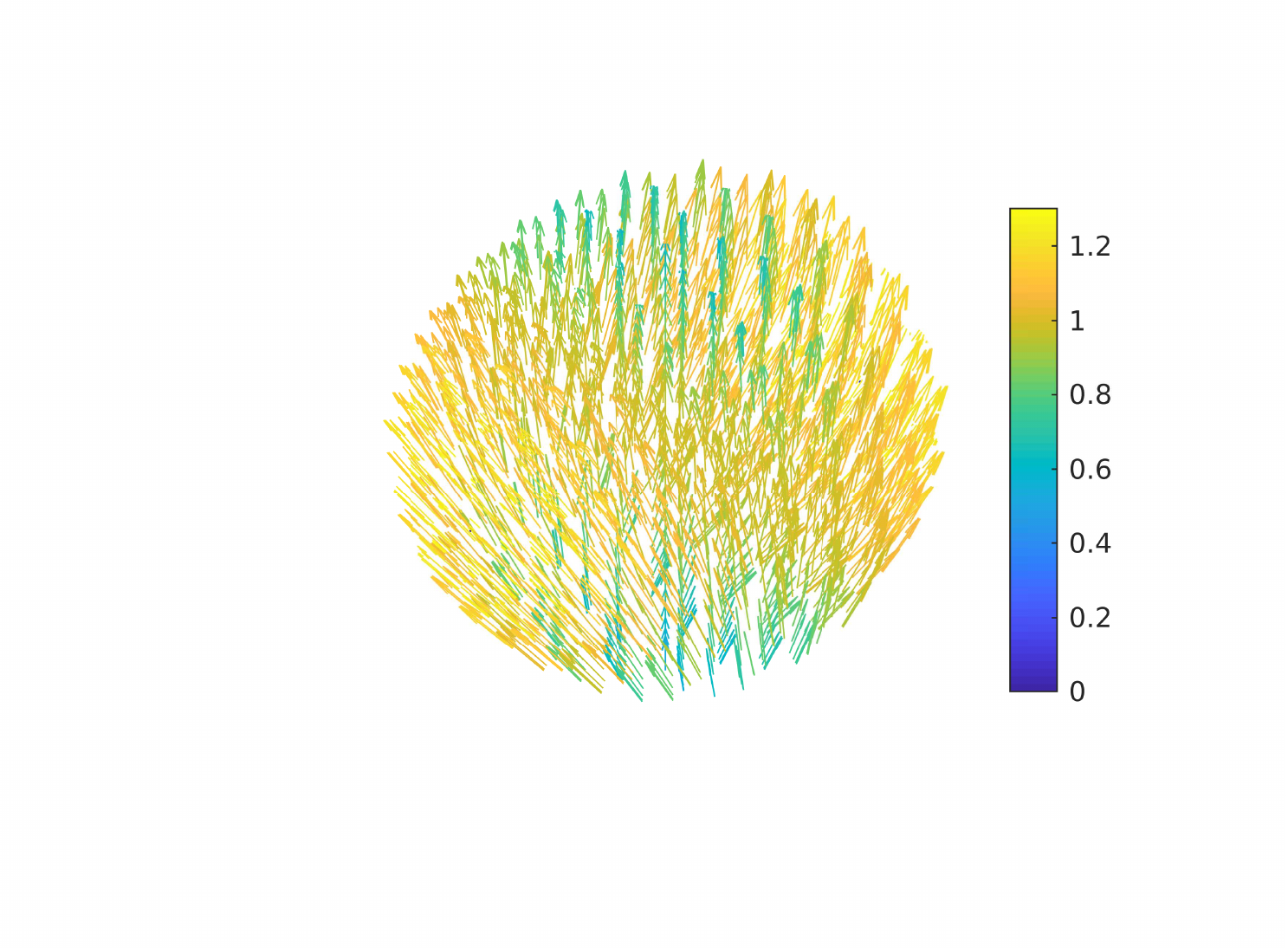}
\end{overpic} 
\end{minipage}
\begin{minipage}{.49\textwidth} 
\begin{overpic}[width=\textwidth,trim={120 70 50 50},clip]{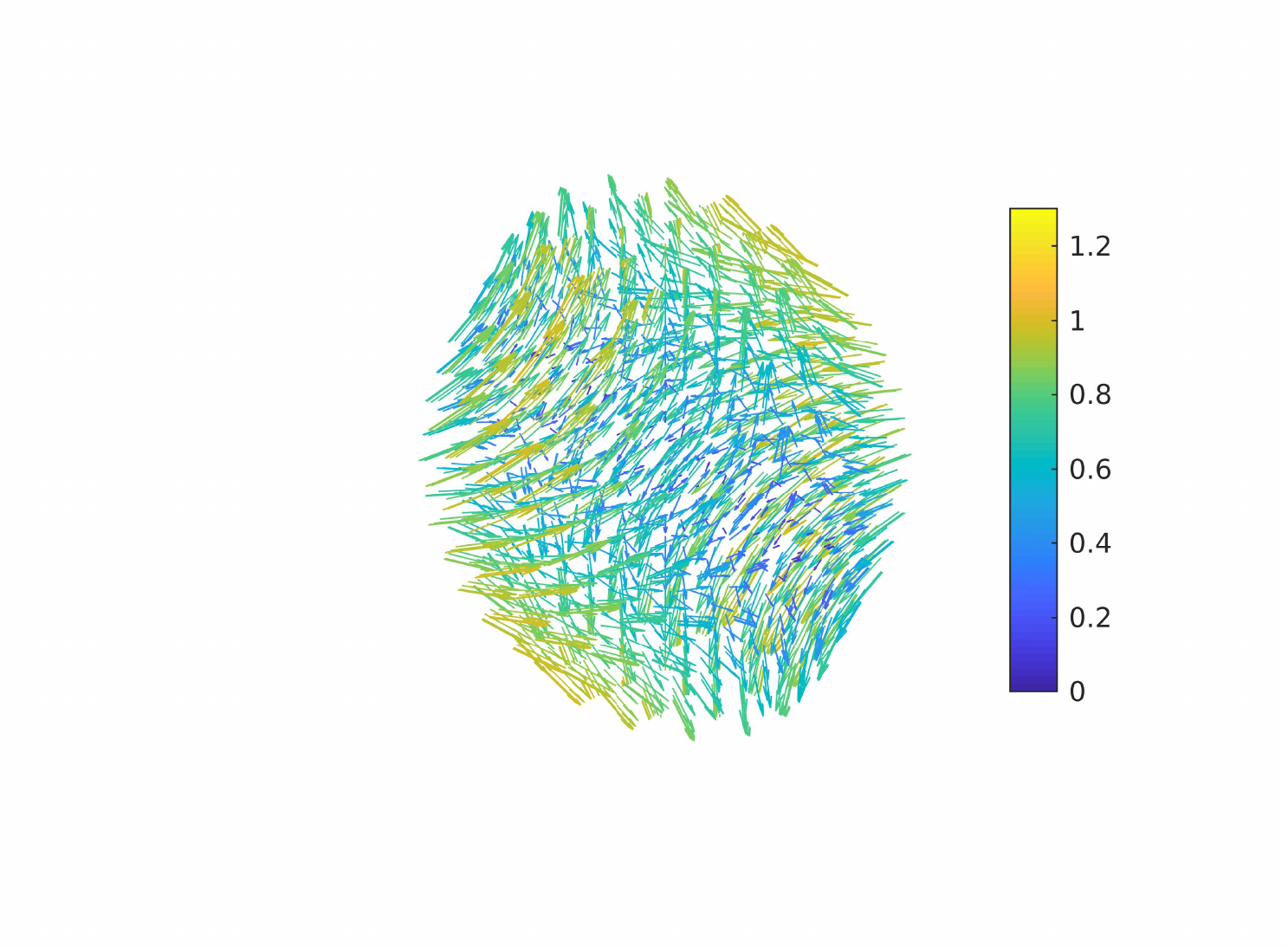}
\end{overpic}
\end{minipage}
\label{fig_Vector_calculus}
\caption{The vector-valued function $V(x,y,z)=(\sin(x),xz,\cos(z))$ (left) and its curl (right), plotted using the \texttt{quiver} command.}
\end{figure}

One can also verify basic vector calculus identities. For example, the divergence theorem asserts that the vector-valued function $v$ satisfies 
\[
\iiint_{B(0,1)}(\nabla\cdot V)\,dV = \oiint_{\partial B} V\cdot\vec{n}\,dS,
\]
where $\hat{n}$ denotes the unit normal vector to $B(0,1)$.  This theorem can be verified by executing the following code: 
\begin{verbatim}
  lhs = sum3(divergence(V));  % Compute volume-integral of div( V )
  nhat = spherefunv.unormal;  % Unit normal vector to surface of ball
  Vbc = V(1,:,:, 'spherical');% Restrict V to the bdy of ball
  rhs = sum2(dot(Vbc, nhat)); % Compute dot-product & surface-integral
\end{verbatim}
The error is $2.2204\times 10^{-15}$.

\subsection{Poloidal-toroidal decomposition}\label{sec_PT_decomposition}
The poloidal-toroidal (PT) decomposition of a smooth divergence-free vector field expresses the field as the sum of two orthogonal fields. The PT decomposition is a well-known tool in fluid dynamics~\cite{horn2015toroidal} and magnetohydrodynamics~\cite{benton1983rapid,boronski2007poloidal,boronski2007poloidal2} simulations to analytically impose an incompressibility condition on flows in cylindrical and spherical geometries. In this section, we describe an algorithm for computing the PT decomposition of a smooth vector field in the ball. 


Given a smooth divergence-free vector field, $V$, defined on the ball, the PT decomposition~\cite{backus1986poloidal} writes $V$ as an orthogonal sum of a poloidal and toroidal field, i.e., 
\[V = P + T, \qquad  \int_{B(0,1)}P\cdot T \,dV = 0.\]
Here, there exist two poloidal and toroidal scalar-valued functions $\Phi$ and $\Psi$ such that $P = \nabla\times \nabla \times (r \Phi \mathbf{e}^r)$ and $T = \nabla\times (r\Psi\mathbf{e}^r)$, where $\mathbf{e}^r$ is the unit radial vector. It can be shown that $\Phi$ and $\Psi$ are unique up to the addition of an arbitrary function of $r$~\cite{backus1986poloidal}. 

A vector field $V$ whose components are expressed in the Cartesian coordinate system $(V_x,V_y,V_z)$ can be converted in spherical coordinates $(V_r,V_\lambda,V_\theta)$ using the following identities:
\begin{align*}
V_r &= \sin\theta(\cos\lambda V_x+\sin\lambda V_y)+\cos\theta V_z,\\
V_\lambda &= -\sin\lambda V_x+\cos\lambda V_y,\\
V_\theta &= \cos\theta(\cos\lambda V_x+\sin\lambda V_y)-\sin\theta V_z.
\end{align*}
Then, any poloidal and toroidal scalars for $V=(V_r,V_{\lambda},V_{\theta})$ satisfy the following relations~\cite{backus1986poloidal}:
\begin{align}
\nabla^2_1\Phi&=-rV_r,\label{eq_Laplace_1}\\
\nabla^2_1\Psi&=\frac{1}{\sin\theta}\left[\frac{\partial}{\partial \theta}(V_{\lambda}\sin\theta)-\frac{\partial}{\partial\lambda}V_{\theta}\right],\label{eq_Laplace_2}
\end{align}
where $\nabla^2_1$ stands for the dimensionless Laplacian defined in the spherical coordinate system $(r,\lambda,\theta)$ by
\[\nabla^2_1 = \frac{1}{\sin\theta}\frac{\partial}{\partial \theta}\left(\sin\theta\frac{\partial}{\partial \theta}\right)+\frac{1}{\sin^2\theta}\frac{\partial^2}{\partial \lambda^2}.\]
After multiplying by $\sin^2\theta$ to remove the singularities, \cref{eq_Laplace_1,eq_Laplace_2} become
\begin{align}
\sin\theta\cos\theta\frac{\partial \Phi}{\partial\theta}+\sin^2\theta\frac{\partial^2\Phi}{\partial\theta^2}+\frac{\partial^2\Phi}{\partial\lambda^2}&=-r\sin^2\theta V_r,\label{eq_poloidal_scalar}\\
\sin\theta\cos\theta\frac{\partial \Psi}{\partial\theta}+\sin^2\theta\frac{\partial^2\Psi}{\partial\theta^2}+\frac{\partial^2\Psi}{\partial\lambda^2}&=-\sin\theta\left[\partial_{\theta}(V_{\lambda}\sin\theta)-\partial_{\lambda}V_{\theta}\right].\label{eq_toroidal_scalar}
\end{align}
Moreover, any smooth function $u$ on the unit ball has a unique interpolant $\tilde{u}$ (\cref{sec_intro}) of the form
\[\tilde{u}(r,\lambda,\theta)\approx\sum_{i=0}^{n-1}\sum_{j=-n/2}^{n/2-1}\sum_{k=-n/2}^{n/2-1} \alpha_{ijk} T_{i}(r)e^{\mathbf{i}j\lambda}e^{\mathbf{i}k\theta}.\]
Thus, \cref{eq_poloidal_scalar,eq_toroidal_scalar} decouple in $\lambda$ and $r$ with this basis. However these equations are not well defined since $P$ and $T$ are unique up to addition of arbitrary functions of $r$. Then, \cref{eq_poloidal_scalar,eq_toroidal_scalar} are solved numerically with the condition that the zero-th Fourier mode in $\lambda$ and $\theta$ is equal to zero. This is equivalent to $\alpha_{i00}=0$ for all $0\leq i\leq n$.

Let $M_{\sin\cos}$ and $M_{\sin^2}$ be the multiplication matrices for $\sin\theta\cos\theta$ and $\sin^2\theta$ in the Fourier basis, $D_n$ the matrix of differentiation with respect to $\theta$, $F$ the tensor of CFF coefficients of $-r\sin^2\theta v_r$ and $G$ the tensor of CFF coefficients of $-\sin\theta[\partial_{\theta}(v_{\lambda}\sin\theta)-\partial_{\lambda}v_{\theta}]$. For example
\[D_n = \diag\left(\left[-\frac{n}{2}i,\cdots,-i,0,i,\cdots,\frac{n}{2}i\right]\right).\]
Equations~\cref{eq_poloidal_scalar,eq_toroidal_scalar} are discretized into
\begin{align*}
\left(M_{\sin\cos}D_n+M_{\sin^2}-j^2I\right)P(i,j,:)&=F(i,j,:)\\
\left(M_{\sin\cos}D_n+M_{\sin^2}-j^2I\right)T(i,j,:)&=G(i,j,:)
\end{align*}
for $0\leq i\leq n-1$ and $-n/2\leq j\leq n/2-1$. $P(i,j,:)$ denotes the vector of Fourier coefficients in $\theta$.
These equations are of the form
\[AX(i,j,:) = B(i,j,:),\]
where $A$ is a sparse banded matrix with bandwidth $m=4$, which is the Fourier number of $\sin\theta\cos\theta$ and $\sin^2\theta$. Then, by Band Gaussian Elimination~\cite{golub2012matrix}, it can be solved in $\mathcal{O}(m^2n)$ operations. Thus, the complexity of the poloidal-toroidal algorithm is $\mathcal{O}(n^3)$, which is linear in the number of interpolation points over the ball. 

As an example we consider the induction equation~\cite[Chap.~2]{davidson2002introduction}, which is one of the equations arising in magnetohydrodynamics and resulting from Maxwell's equations,
\begin{align}
\frac{\partial B}{\partial t}&=\nabla\times(u\times B)+D\nabla^2 B, \label{eq_maxwell_1} \\
\nabla\cdot B&=0, \label{eq_maxwell_2}
\end{align}
where $B$ denotes the magnetic field, $u= \nabla\times[e^{-5(x^2+y^2+z^2)}(x^2,y^2,xz)]$ is the velocity of particles, and $D=1/3000$ is the diffusion constant. We use the poloidal-toroidal decomposition to ensure that the divergence-free condition on $B$ is satisfied and decouple \cref{eq_maxwell_1} into the following equations on the poloidal and toroidal scalars $\Phi_B$ and $\Psi_B$, 
\begin{align}
\frac{\partial \Phi_B}{\partial t}&=\Phi_{\nabla\times(u\times B)}+D\nabla^2 \Phi_B, \label{eq_maxwell_3}\\
\frac{\partial \Psi_B}{\partial t}&=\Psi_{\nabla\times(u\times B)}+D\nabla^2 \Psi_B. \label{eq_maxwell_4}
\end{align}
Then, according to the IMEX-BDF1 time-stepping  scheme (see~\cref{sec_ex_helmholtz}), at each time step we compute the poloidal-toroidal decomposition of the nonlinear term $\nabla\times(u\times B)$ and solve two Helmholtz equations:
\begin{align*}
\nabla^2\Phi_B^{n+1}+K^2\Phi_B^{n+1}&=K^2\Phi_B^n+\frac{1}{D}\Phi_{\nabla\times(u\times B)}^n,\\
\nabla^2\Psi_B^{n+1}+K^2\Psi_B^{n+1}&=K^2\Psi_B^n+\frac{1}{D}\Psi_{\nabla\times(u\times B)}^n,
\end{align*}
where $\Phi_B^n$ (resp. $\Psi_B^n$) denotes the magnetic poloidal (resp. toroidal) scalar at time $t = n\Delta t$, $\Delta t = 5\times 10^{-2}$ is the time step, and $K^2 = -1/(D\Delta t)$. We choose the initial magnetic field $B= \nabla\times[ze^{-5(x^2+y^2+z^2)}(x,y,z)]$ and impose homogeneous Dirichlet boundary conditions on the poloidal and toroidal scalars, which are computed at each time step using the following code snippet:
\begin{verbatim}
  B = ballfunv.PT2ballfunv(Phi_B, Psi_B);% Compute B from P and T
  N = curl(cross(u, B));                 % Nonlinear term
  [Phi_N, Psi_N] = PTdecomposition(N);   % PT decomposition of N
  % Solve the  toroidal and poloidal equation
  Phi_B = helmholtz(K^2*Phi_B+Phi_N/D, K, @(x,y,z)0, 100);
  Psi_B = helmholtz(K^2*Psi_B+Psi_N/D, K, @(x,y,z)0, 100);
\end{verbatim}

\begin{figure}[htbp]
\begin{minipage}{.32\textwidth} 
\begin{overpic}[width=\textwidth,trim={120 70 50 50},clip]{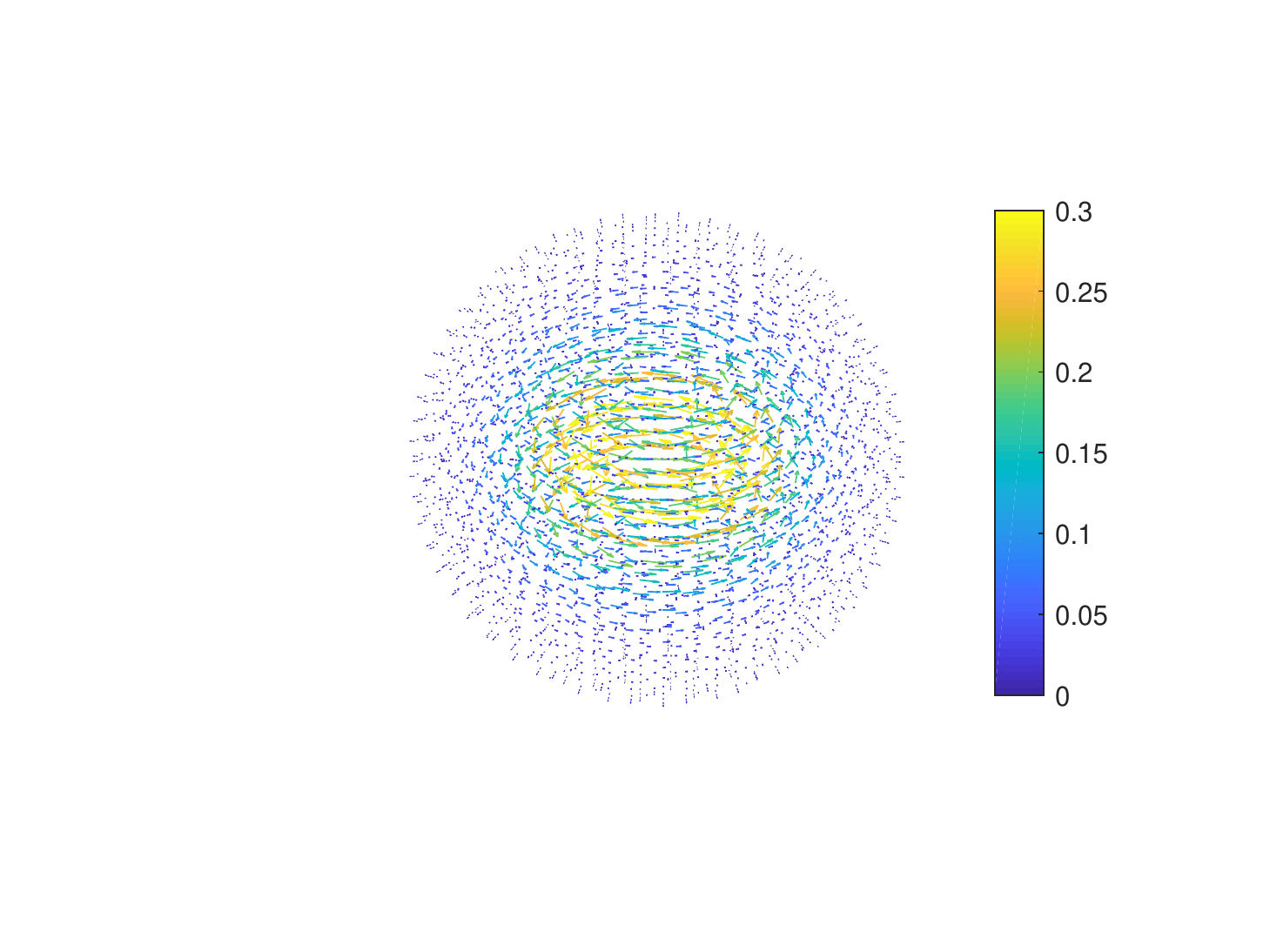}
\put(0,65){(a)}
\end{overpic} 
\end{minipage}
\begin{minipage}{.32\textwidth} 
\begin{overpic}[width=\textwidth,trim={120 70 50 50},clip]{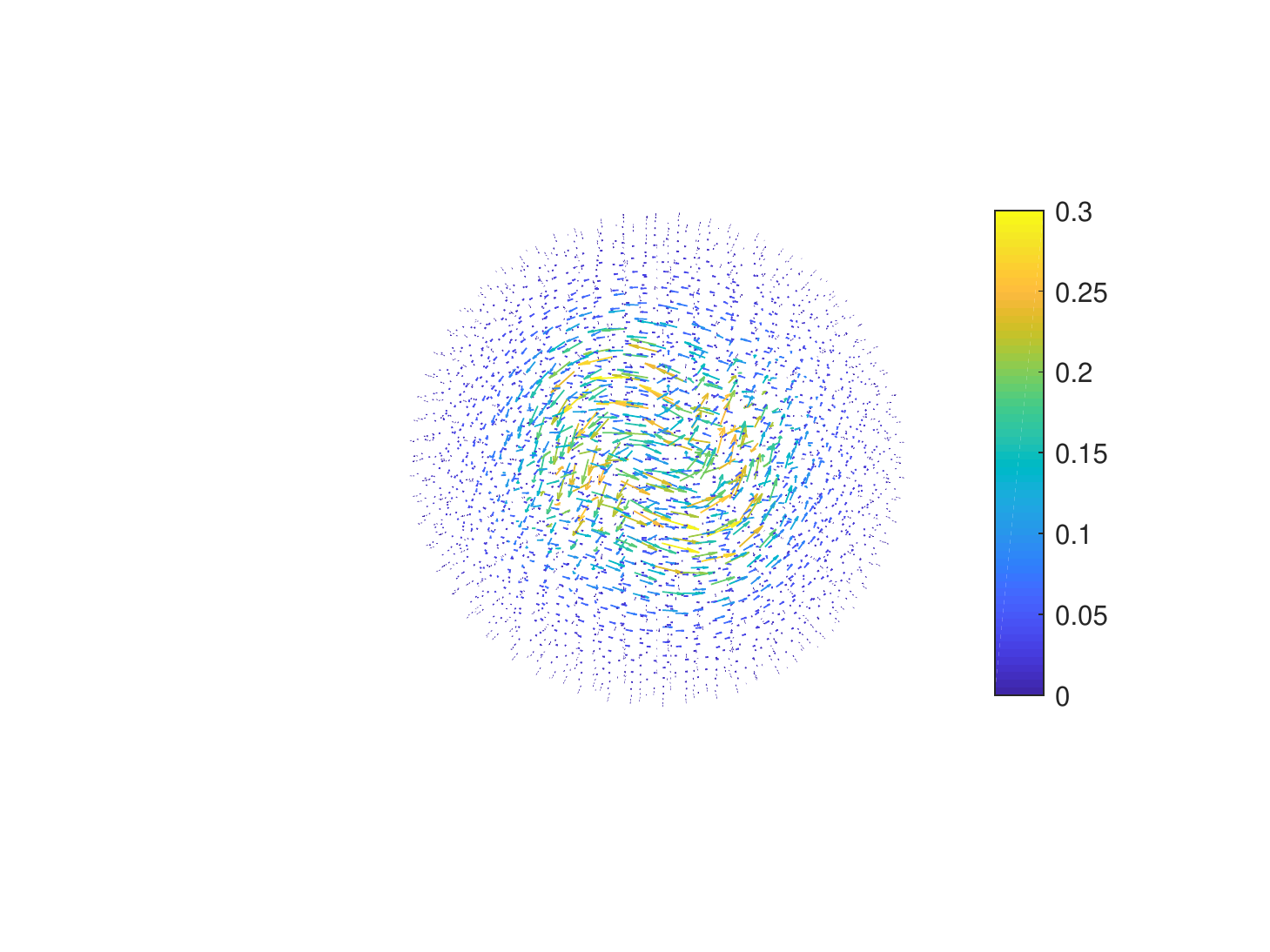}
\put(0,65){(b)}
\end{overpic} 
\end{minipage}
\begin{minipage}{.32\textwidth} 
\begin{overpic}[width=\textwidth,trim={120 70 50 50},clip]{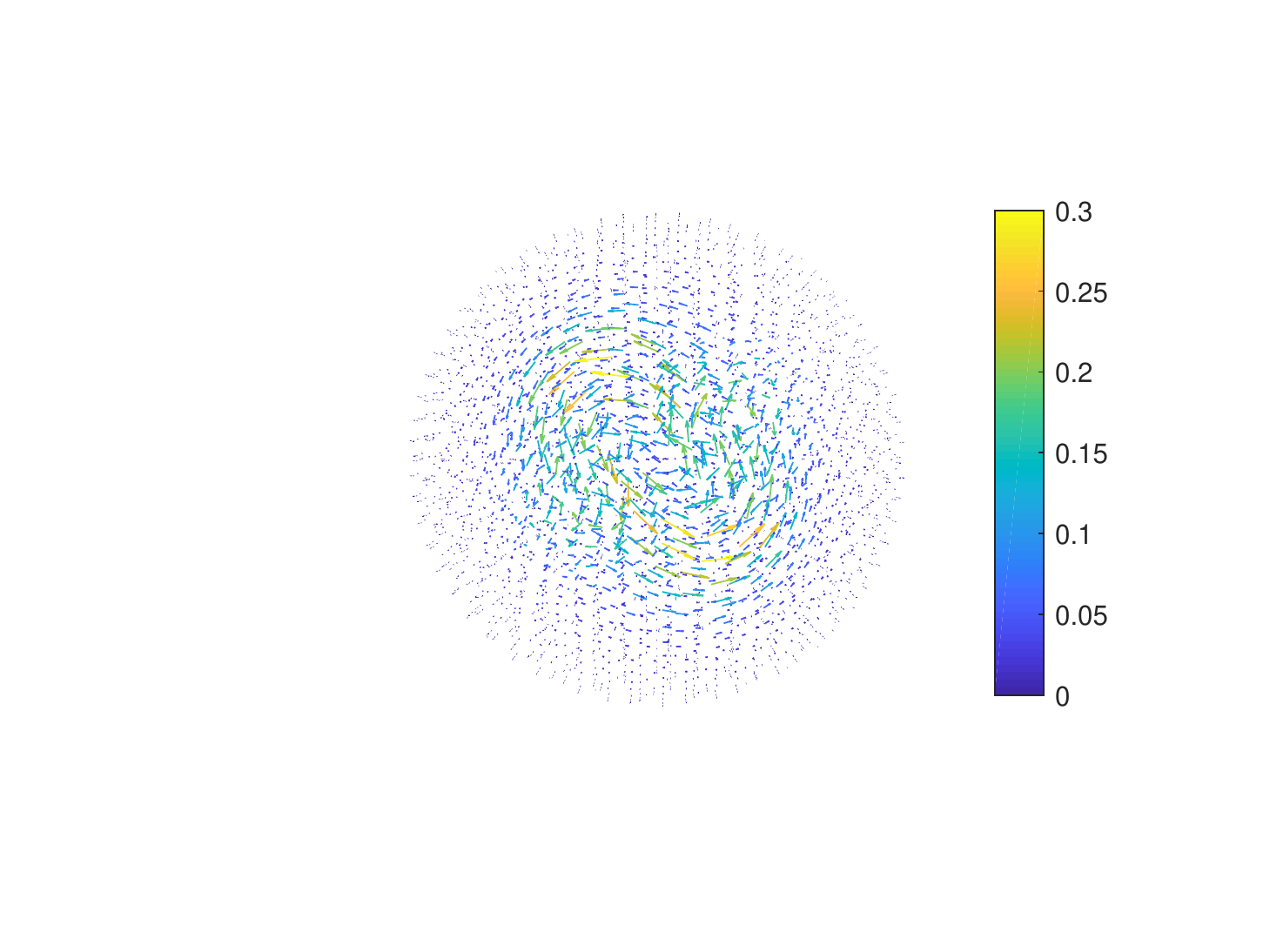}
\put(0,65){(c)}
\end{overpic} 
\end{minipage}
\label{fig_maxwell}
\caption{Solution to the induction equation \eqref{eq_maxwell_1}-\eqref{eq_maxwell_2} at time $t = 0$ (see (a)), $t = 1$ (see (b)), and $t = 2$ (see (c)).}
\end{figure}

The magnetic field $B$ is illustrated at different time steps in \cref{fig_maxwell}.

\subsection{Helmholtz--Hodge decomposition}
The Helmholtz--Hodge decomposition (HHD) has been an important tool in computational fluid dynamics since the introduction of projection methods by Chorin~\cite{chorin1967numerical,chorin1968numerical,chorin1969convergence} to solve the Navier--Stokes equations for incompressible fluids. The decomposition is then used to preserve the divergence-free property of the velocity field during the computation of the solution. Applications of the HHD also arise in computer graphics and visualization of incompressible fluids such as water~\cite{bridson2015fluid,pharr2016physically,tan2009physically}. The HHD has also been exploited in the field of computer vision and robotics to analyze cardiac videos~\cite{guo2006cardiac} and find singularities in fingerprints images~\cite{gao2010singular}. A survey of applications is available in~\cite{bhatia2013helmholtz}. Ballfun has a \texttt{HelmholtzDecomposition} command, which computes and returns the Helmholtz--Hodge decomposition of a smooth vector field using the algorithm described in this section.

The Helmholtz--Hodge decomposition~\cite{bhatia2013helmholtz} states that any smooth vector field $v$ on the unit ball can be decomposed into a sum of a solenoidal and irrotational fields
\begin{equation}
\label{eq_Helmholtz}
V = \nabla f + \psi,
\end{equation}
where $\psi$ is a divergence-free vector field. Moreover, this decomposition can be made unique by imposing that the incompressible component, $\psi$, is tangent to the boundary.  Moreover, this condition is equivalent to a Neumann boundary condition on the scalar function $f$. That is $\vec{n}\cdot\nabla f=\vec{n}\cdot V$~\cite{bhatia2013helmholtz}. The first step of the algorithm implemented in Ballfun consists of taking the divergence of $V$ in \cref{eq_Helmholtz} to obtain the Poisson equation
\[\nabla^2 f = \nabla\cdot V\]
with the Neumann boundary conditions given by
\[\vec{n}\cdot \nabla f = \left.\frac{\partial f}{\partial r}\right|_{\partial B} = \vec{n}\cdot V.\]
We then obtain the incompressible component using the following equality:
\[\psi:=V-\nabla f.\]
Finally, the poloidal and toroidal scalars of $\psi$ are computed using the algorithm described in \cref{sec_PT_decomposition} and the command \texttt{HelmholtzDecomposition}
returns the scalar function $f$ together with the poloidal and toroidal scalars of $\psi$.

\begin{figure}[htbp]
\centering
\begin{minipage}{.32\textwidth} 
\begin{overpic}[width=\textwidth,trim={120 70 50 50},clip]{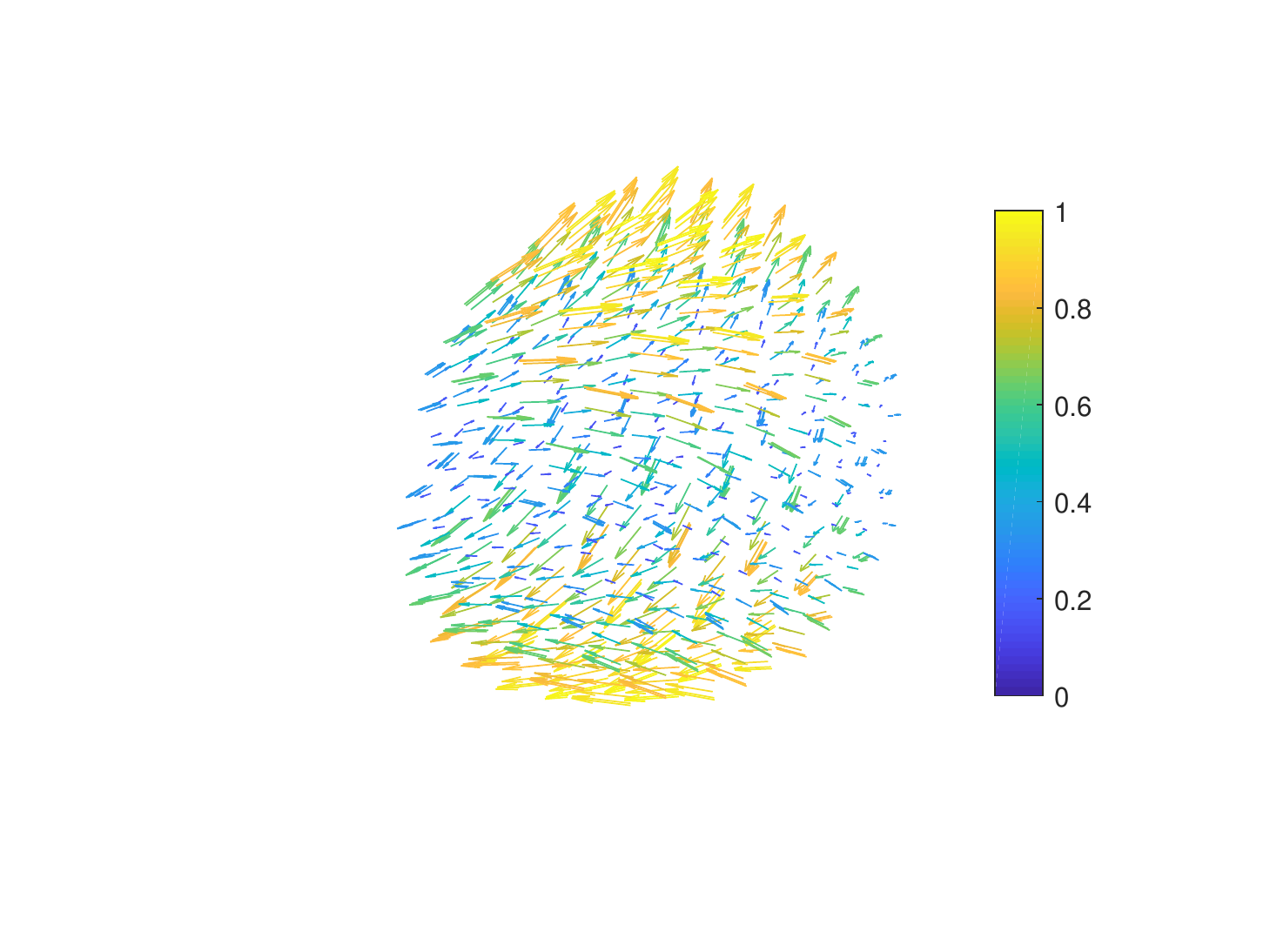}
\put(0,65){(a)}
\end{overpic} 
\end{minipage}
\begin{minipage}{.32\textwidth} 
\begin{overpic}[width=\textwidth,trim={120 70 50 50},clip]{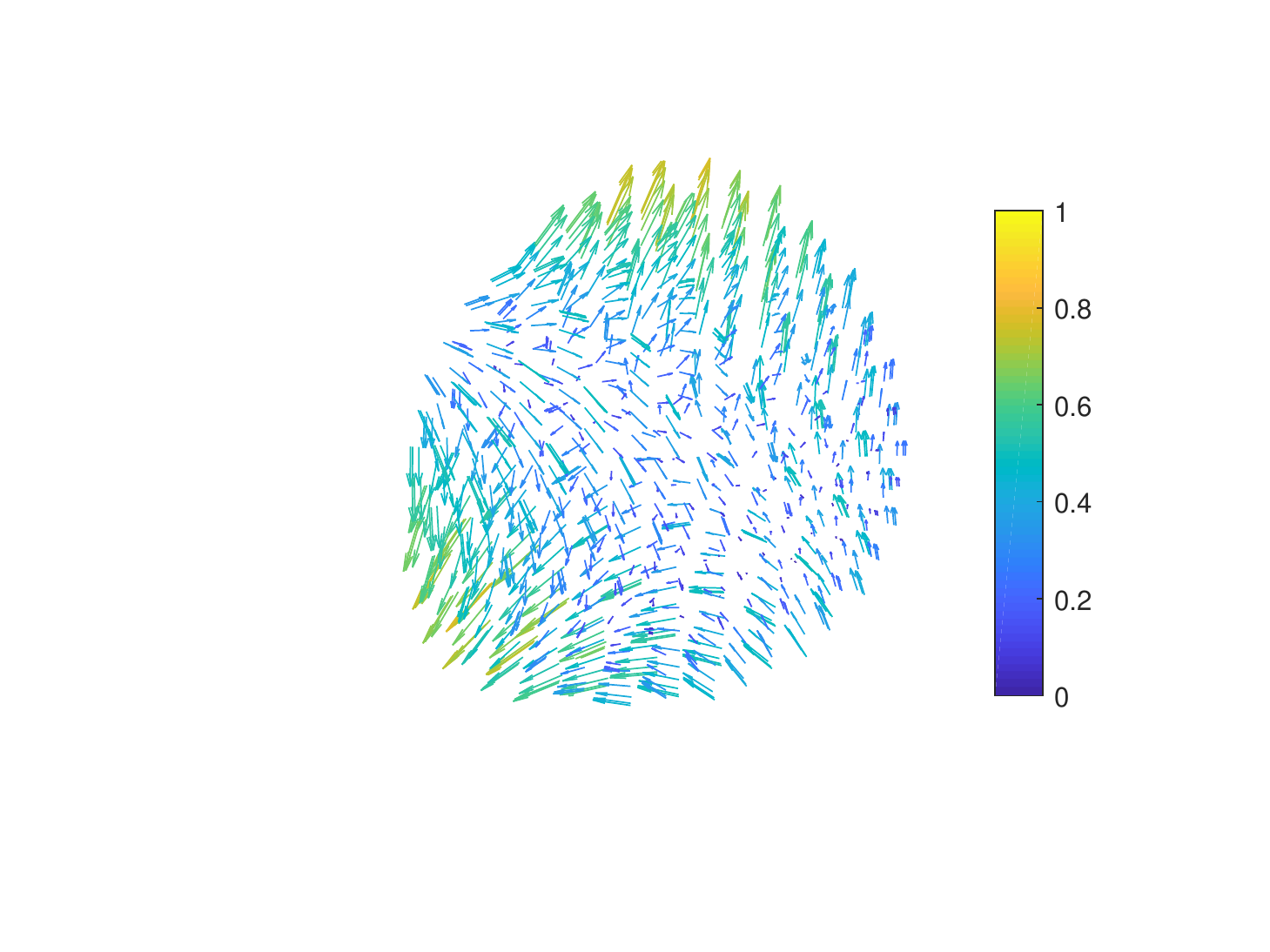}
\put(0,65){(b)}
\end{overpic} 
\end{minipage}
\begin{minipage}{.32\textwidth} 
\begin{overpic}[width=\textwidth,trim={120 70 50 50},clip]{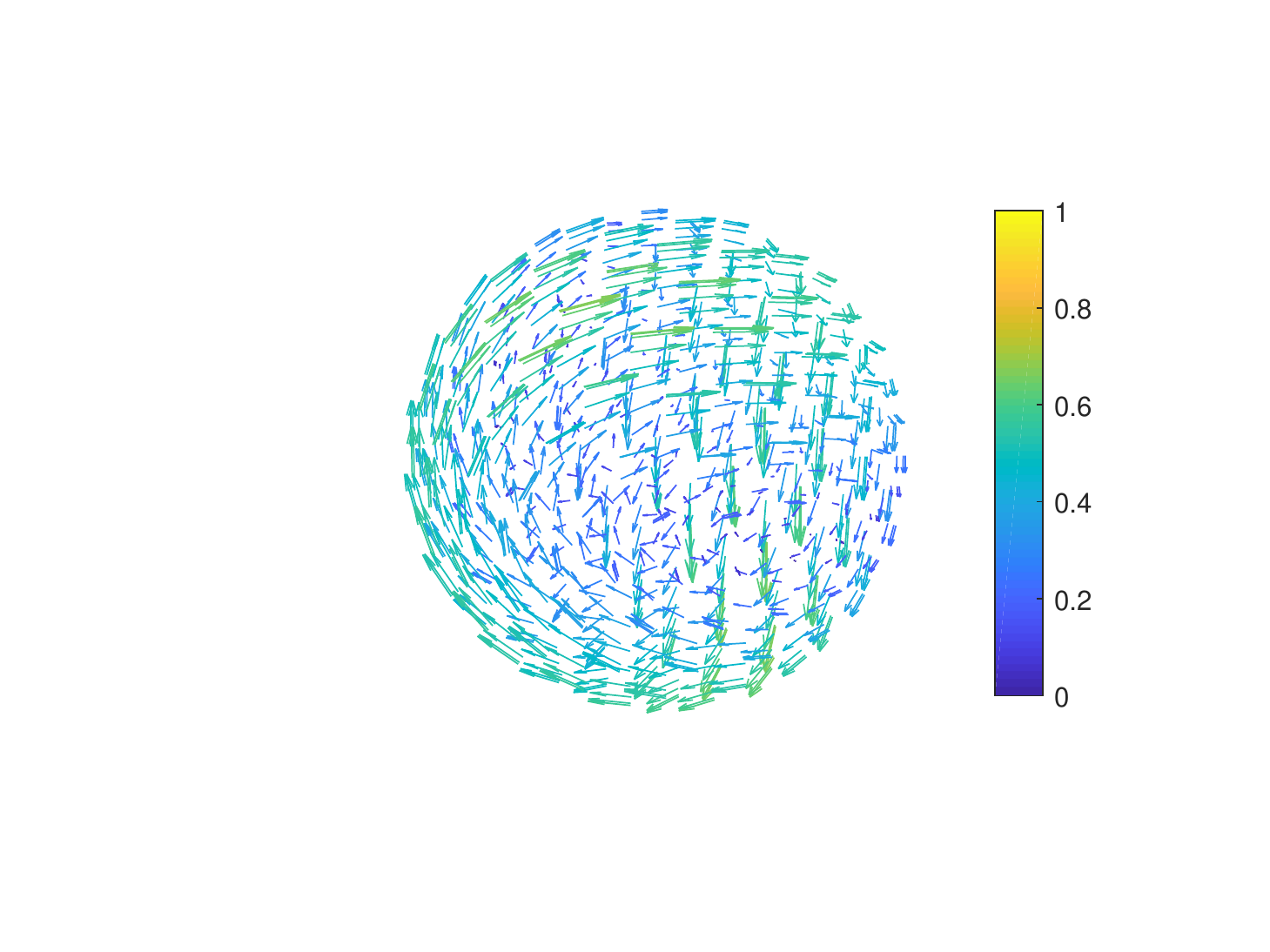}
\put(0,65){(c)}
\end{overpic} 
\end{minipage}
\label{fig_HHD}
\caption{The vector field $V(x,y,z) = (\cos(xy)z,\sin(xz),yz)$ (see (a)), together with its Helmholtz--Hodge decomposition $\nabla f$ (see (b)) and $\Psi$ (see (c)). The decomposition is computed using \texttt{HelmholtzDecomposition} and plotted with the \texttt{quiver} command.}
\end{figure}

\cref{fig_HHD} shows the Helmholtz--Hodge decomposition of the vector-valued function $V(x,y,z) = (\cos(xy)z,\sin(xz),yz)$ computed by Ballfun using the following code:
\begin{verbatim}
  % Define the vector field v
  v = ballfunv(@(x,y,z)cos(x.*y).*z, @(x,y,z)sin(x.*z), @(x,y,z)y.*z);
  % Compute the Helmholtz-Hodge decomposition of v
  [f, Ppsi, Tpsi] = HelmholtzDecomposition( v );
  % Recover Psi from it poloidal and toroidal scalars
  Psi = ballfunv.PT2ballfunv(Ppsi,Tpsi);
\end{verbatim}

\section{Conclusions}
The analogue of the DFS method for the ball is exploited to impose BMC structure on functions and represent them by Chebyshev--Fourier--Fourier series. A collection of fast and spectrally accurate algorithms is developed for differentiation, rotation, solving the Helmholtz equation, vector calculus, poloidal-toroidal decomposition, and Helmholtz--Hodge decomposition. These ideas have been implemented in Ballfun, which is part of the freely available software Chebfun.

\section*{Acknowledgments}
We thank Vassilios Dallas and Jonasz S{\l}omka for discussions on the poloidal-toroidal and Helmholtz--Hodge decomposition. We also thank Nick Trefethen, Heather Wilber, and Grady Wright for comments and suggestions on the Ballfun software and the paper.

\newpage
\bibliographystyle{siamplain}
\bibliography{biblio}
\end{document}